\documentclass{amsart}
\usepackage{amssymb}
\newtheorem{theorem}{Theorem}[section]

\newtheorem{remark}[theorem]{Remark}
\newtheorem{lemma}[theorem]{Lemma}
\newtheorem{prop}[theorem]{Proposition}
\newtheorem{corol}[theorem]{Corollary}
\begin{document}
\renewcommand{\v}{\varepsilon} \newcommand{\p}{\rho}
\newcommand{\m}{\mu}
\def\im{{\bf im}}
\def\re{{\bf re}}
\def\e{{\bf e}}
\def\a{\alpha}
\def\ve{\varepsilon}
\def\b{\beta}
\def\D{\Delta}
\def\d{\delta}
\def\f{{\varphi}}
\def\ga{{\gamma}}
\def\l{\lambda}
\def\L{\Lambda}
\def\lo{{\bf l}}
\def\s{{\bf s}}
\def\A{{\bf A}}
\def\B{{\bf B}}
\def\cB{{\mathcal {B}}}
\def\C{{\mathbb C}}
\def\F{{\bf F}}
\def\G{{\mathfrak {G}}}
\def\g{{\mathfrak {g}}}
\def\H{{\mathcal {H}}}
\def\O{\Omega}
\def\M{{\mathcal {M}}}
\def\N{{\mathcal {N}}}
\def\R{{\mathcal {R}}}
\def\U{{\mathcal {U}}}
\def\Z{{\mathbb Z}}
\def\P{{\mathcal {P}}}
\def\GVM{ GVM }
\def\iff{ if and only if  }
\def\add{{\rm add}}
\def\ld{\ldots}
\def\vd{\vdots}
\def\sl{{\rm sl}}
\def\mod{{\rm mod}}
\def\len{{\rm len}}
\def\cd{\cdot}
\def\dd{\ddots}
\def\q{\quad}
\def\qq{\qquad}
\def\ol{\overline}
\def\tl{\tilde}
\def\nn{\nonumber}
\title[Integrable modules]
{Integrable  modules for affine Lie
superalgebras}
\author{Senapathi Eswara Rao}
\address{Tata Institute of Fundamental Research, Mumbai, India}
\email{senapati@math.tifr.res.in}
\author{Vyacheslav Futorny}
\address{
Institute of Mathematics, University of S\~ao Paulo, Caixa Postal 66281
CEP 05315-970,
S\~ao Paulo, Brazil}
\email{futorny@ime.usp.br}
\thanks{The second author  is supported in part by the CNPq grant 307812/2004-9 and by the Fapesp grant 2005/60337-2.}

\date{}

\begin{abstract}
 Irreducible nonzero level modules  with finite-dimensional weight spaces are
 discussed
for non-twisted affine Lie superalgebras. A complete classification of such modules is obtained for superalgebras of
type $A(m,n)^{\hat{}}$ and $C(n)^{\hat{}}$ using Mathieu's
classification of cuspidal modules over simple Lie algebras. In
other cases  the classification problem is reduced to the
classification of  cuspidal modules over  finite-dimensional
cuspidal Lie superalgebras described by Dimitrov, Mathieu and
Penkov. Based on these results a complete classification of
irreducible integrable (in the sense of Kac and Wakimoto) modules
is obtained by showing that any such module is  highest weight in
which case the problem was solved by  Kac and Wakimoto.
\end{abstract}
\maketitle
\section{ Introduction}
Let $\g=\g_0\oplus \g_1$ be a Lie superalgebra over $\C$, i.e. 
\begin{itemize}
\item $[\g_{\ve}, \g_{\ve'}]\subset \g_{\ve+\ve'(\mod 2)}$;
\item $[x,y]=-(-1)^{|x||y|}[y,x]$;
\item $[x,[y,z]]=[[x,y],z]+(-1)^{|x||y|}[y,[x,z]]$,
\end{itemize}

where $|x|=\ve$ if $x\in \g_{\ve}$.
We will assume that $\g$ is finite-dimensional basic classical Lie superalgebra, i.e. a simple Lie superalgebra with a non-degenerate even supersymmetric invariant bilinear form and with reductive $\g_0$. The classification of such algebras was obtained in \cite{K}.

Let $V=V_0\oplus V_1$ be a $\Z_2$-graded vector space, $\dim
V_0=m$, $\dim V_1=n$. Then the endomorphism algebra $End V$ is an associative superalgebra with a
natural $\Z_2$-gradation. 
Defining the Lie bracket as $[A,B]=AB-(-1)^{|A||B|}BA$ we
make $End V$ into a Lie superalgebra $\mathfrak{gl}(m,n)$. The {\em supertrace} is a linear
function $str:\mathfrak{gl}(m,n)\rightarrow \C$ such that $str
id_V=m-n$ and $str [A,B]=0$, $A,B\in \mathfrak{gl}(m,n)$. Denote
$$\mathfrak{sl}(m+1,n+1)=\{A\in \mathfrak{gl}(m+1,n+1)|str A=0\}.$$ If $m\neq n$ this
is the superalgebra of type $A(m,n)$. In case $m=n$,
$\mathfrak{sl}(n+1,n+1)$ has a one-dimensional ideal consisting of
scalar matrices. Its quotient $\mathfrak{psl}(n+1,n+1)$ is basic classical Lie
superalgebra of type $A(n,n)$. 
Let
$F=\Big({\begin{array}{c} I_m \quad 0\\
0 \quad -I_n \end{array}}\Big)$, $C=\Big({\begin{array}{c} C_1 \quad 0\\
0 \quad C_2 \end{array}}\Big)$, $C_1^t=C_1,$ $C_2^t=-C_2$.
Then 
$$\mathfrak{osp}(m,2n)_a=\{A\in
\mathfrak{gl}(m,n)_a|F^aA^tC+CA=0\},\, a=0,1.$$

Let $$B(m,n)=\mathfrak{osp}(2m+1,2n),$$
$$C(n)=\mathfrak{osp}(2,2n),$$ $$D(m,n)=\mathfrak{osp}(2m,2n), m>1.$$

The series $A(m,n)$, $B(m,n)$, $D(m,n)$, $C(n)$ and 
special Lie superalgebras
 $D(2,1;a)$, $G(3)$, $F(4)$ (cf. Section~5) form a complete list of basic classical Lie superalgebras.

For the purposes of this paper  we will consider the Lie superalgebra 
$\mathfrak{pgl}(n+1,n+1)$ instead of $\mathfrak{psl}(n+1,n+1)$ as a superalgebra of type $A(n,n)$.

Let $\G=\G_0\oplus \G_1$ be the corresponding non-twisted affine Lie superalgebra, i.e. $1$-dimensional central extension of the loop superalgebras $\g\otimes \C[t,t^{-1}]$ with a degree derivation $d$, $d(x\otimes t^n)=n(x\otimes t^n)$. Denote by
$\H$  a Cartan subalgebra of $\G_0$ and let $c \in \H$ be  the central element of
$\G$. Also let
$\D$ be
the root system of $\G$ with respect to $\H$, $\D^{\re}$ (respectively $\D^{\im}$)  the set
of real (respectively imaginary) roots, $\pi=\pi_0\cup \pi_1$  a basis of $\D$
and $Q$  the free abelian group generated by $\pi$. 

If $V$ is a $\G$-module and $\l\in \H^*$ we set $V_\l=\{v\in V |\  hv=\l(h)v, {\rm \ for\ all\ } h\in \H\}$.
If $V_\l$ is non-zero then $\l$ is a {\it weight} of $V$ and
$V_\l$ is the corresponding weight space.  We denote by $P(V)$
   the set of all weights of $V$.  A $\G$-module is  {\it weight}
if $V=\bigoplus V_\l$, $\l\in P(V)$.

Similarly one defines weight modules for finite-dimensional reductive Lie algebras and Lie superalgebras.
  A weight module for a reductive Lie algebra $\cB$ is
 {\it cuspidal}
if  all its  root elements act injectively. In this case all simple components  of $\cB$ are of type $A$ and $C$
\cite{Fe}. Irreducible cuspidal modules with finite-dimensional weight spaces were classified by Mathieu \cite{M}.
   The concept of cuspidality has been extended to weight irreducible modules over finite-dimensional
 Lie superalgebras by Dimitrov, Mathieu and Penkov in \cite{DMP}.   Such module $V$ is
cuspidal if the monoid generated by the roots with the injective
action of the corresponding root elements, is a subgroup of finite index in $Q$, which is equivalent to the fact
that $V$ is not parabolically induced (cf. \cite{DMP}, Corollary 3.7).

If $c$ acts as a scalar on a $\G$-module $V$ then this scalar is
called the level of $V$. Denote by $K(\G)$ the category of weight
$\G$-modules of nonzero levels with finite-dimensional weight
spaces. The first goal of the present paper is to describe the
irreducible modules in the category $K(\G)$. For an affine Lie
algebra $\cB$ the classification of irreducible modules in
$K(\cB)$ was obtained in \cite{FT}. It was shown that any such
module is a quotient of a
 module induced from an irreducible cuspidal module over a finite-dimensional reductive Lie subalgebra.

We show that in the case of affine Lie superalgebra $\G$ the classification of
  irreducible $\G$-modules in $K(\G)$ is reduced to the classification  of irreducible cuspidal modules over
  {\it cuspidal} Levi subsuperalgebras. The cuspidal  Levi subsuperalgebras of affine Lie superalgebras are
essentially the same as of finite-dimensional classical Lie superalgebras in which case
they where described in  \cite{DMP} and \cite{DMP1}.

A subset $\P\subset \D$ is called a {\em parabolic subset} if
$\P$ is additively closed and $\P\cup -\P=\D$. This concept
in the affine Lie algebras setting was introduced in \cite{F1}.
For a parabolic subset $\P$ denote $\P^{\pm}=\pm (\P\setminus
(-\P))$ and $\P^0= \P\cap (-\P)$. This induces 
  the corresponding
 triangular decomposition of $\G$:
 $\G=\G_{\P}^-\oplus (\G_{\P}^0+\H)\oplus \G_{\P}^+$. Denote $\G_{\P}=(\G_{\P}^0+\H)\oplus \G_{\P}^+$.

Now we state our main result.

\medskip

 \noindent {\bf Theorem A}. {\it Let $V$ be an irreducible module in $K(\G)$. Then there exists
a parabolic subset $\P\subset \D$
  and an irreducible weight module $N$ over the corresponding parabolic subalgebra $\G_{\P}$
such that $V$ is parabolically induced from $\G_{\P}$. Moreover, $N$ is irreducible
cuspidal module over the Levi subsuperalgebra of $\G_{\P}$}.

 \medskip

  If $\G$ is of type $A(m,n)^{\hat{}}$ or $C(n)^{\hat{}}$ then combining Theorem A with Mathieu's
classification of irreducible
 cuspidal modules for Lie algebras of type $A$ and $C$
we obtain a classification of irreducible modules in $K(\G)$ (up to the Weyl group action).
   In all other cases the classification is reduced to a finite indeterminacy by
   \cite{DMP},Proposition~6.3. We should point that there might be different parabolic subalgebras that give the same irreducible module.  It would be interesting to study this question in order to get a complete classification.

A parabolic subset $\P$ is {\em standard} (with respect to $\pi$)
if $\P^+\subset \D_+$ and $\P^0$ is the closure of $S$ in $\D$ for some $S\subset \pi$. In
this case we write $\P=\P_S$ and $\G_{\P}^0=\G(S)$. 
In Section~6 we prove the following  stronger version of Theorem~A.

\medskip

\noindent {\bf Theorem B}. {\it
Modules $L_{\P}(N)$ with irreducible cuspidal $\G_{\P}^0$-module
$N$ exhaust all irreducible objects in $K(\G)$. Moreover $\P$ can
be chosen in the following way:
 \begin{itemize}
 \item[(1)] If $\G$ is of type $A(m,n)^{\hat{}}$  then $\P=\P_S$, $S$ is any  subset of $\pi_0$,
 $\G(S)$ is isomorphic to a finite-dimensional semisimple Lie subalgebra
 of $\mathfrak{sl}(m)\oplus \mathfrak{sl}(n)$.
  \item[(2)] If $\G$ is of type $C(n)^{\hat{}}$ ($n>2$) then
either $\P=\P_S$, $S$ is any  subset of $\pi_0$ and $\G(S)$ is
isomorphic to a finite-dimensional semisimple Lie subalgebra of
$\mathfrak{sp}(2n-2)$, or $\G_{\P}^0$ is
 isomorphic to a semisimple subalgebra of
  $\mathfrak{sp}(2k-2)\oplus \mathfrak{sp}(2m-2)$ with
$k+m=n$.
 \item[(3)] If $\G$ is of type $B(m,n)^{\hat{}}$ ($m>0$)
then $\G_{\P}^0$ is isomorphic to one of the following
(super)algebras or their direct summands:
$$\mathfrak{A}_n\oplus \mathfrak{B}_m; \mathfrak{osp}(4,2n)\oplus
\mathfrak{B}_{m-2};
 \mathfrak{osp}(6,2n)\oplus \mathfrak{B}_{m-3};$$
 $$\mathfrak{A}_k\oplus \mathfrak{osp}(2m+1,2(n-k)), k=0, \ldots, n-1, m=1,2;$$
 $$\mathfrak{sp}(2i)\oplus \mathfrak{sp}(2n-2i)\oplus \mathfrak{sl}(2)\oplus \mathfrak{sl}(2)
\oplus \mathfrak{B}_{m-2};$$
 $$\mathfrak{sp}(2i)\oplus \mathfrak{sp}(2n-2i)\oplus \mathfrak{sl}(3)
\oplus \mathfrak{B}_{m-3}, \mathfrak{sp}(2i)\oplus
\mathfrak{sp}(2n-2i)
 \oplus \mathfrak{B}_{m},$$
$i=0, \ldots, n-1$,
 where $\mathfrak{A}_k$ is a semisimple
subalgebra of $\mathfrak{sp}(2k)$ and $\mathfrak{B}_k$ is a
semisimple subalgebra of $\mathfrak{so}(2k+1)$ with simple
components of type $A$ and $C_2$.

\item[(4)] If $\G$ is of type $B(0,n)^{\hat{}}$  then $\P=\P_S$ and
 $\G(S)$ is a finite-dimensional semisimple Lie subalgebra of $\mathfrak{sp}(2n)$ if
$S$ does not contain the odd
simple root,
  or
$\G(S)=\mathfrak{A}\oplus \mathfrak{osp}(1,2(n-i))$ for some $i\in
\{0, \ldots, n-1\}$, otherwise, where $\mathfrak{A}$ is a
semisimple subalgebra of $\mathfrak{sp}(2i)$.

 \item[(5)] If $\G$ is of type $D(m,n)^{\hat{}}$  then
$\G_{\P}^0$ is isomorphic to one of the following (super)algebras
or their direct summands:
$$\mathfrak{A}_n\oplus \mathfrak{D}_m; \mathfrak{osp}(4,2n)\oplus
\mathfrak{D}_{m-2};
 \mathfrak{osp}(6,2n)\oplus \mathfrak{D}_{m-3};$$
 $$\mathfrak{A}_k\oplus \mathfrak{osp}(2m,2(n-k)), k=0, \ldots, n-1, m=2,3;$$
 $$\mathfrak{sp}(2i)\oplus \mathfrak{sp}(2n-2i)\oplus \mathfrak{sl}(2)\oplus \mathfrak{sl}(2)
\oplus \mathfrak{D}_{m-2};$$
 $$\mathfrak{sp}(2i)\oplus \mathfrak{sp}(2n-2i)\oplus \mathfrak{sl}(3)
\oplus \mathfrak{D}_{m-3}, \mathfrak{sp}(2i)\oplus
\mathfrak{sp}(2n-2i)
 \oplus \mathfrak{D}_{m},$$
 $ i=0, \ldots, n-1,$ where
$\mathfrak{A}_k$ is a semisimple subalgebra of $\mathfrak{sp}(2k)$
and $\mathfrak{D}_k$ is a semisimple subalgebra of
$\mathfrak{so}(2k)$ with simple components of type $A$.

\item[(6)] If $\G$ is of type $D(2,1;a)^{\hat{}}$ then $\P=\P_S$,
$$\G(S)\in \{\mathfrak{sl}(2), \mathfrak{sl}(2)
\oplus \mathfrak{sl}(2),
  \mathfrak{sl}(2)\oplus \mathfrak{sl}(2)\oplus \mathfrak{sl}(2), D(2,1;a)\}.$$
\item[(7)] If $\G$ is of type $G(3)^{\hat{}}$ then $\G_{\P}^0$ is isomorphic
to one of the following \\(super)algebras:
$$ \{\mathfrak{sl}(2), \mathfrak{sl}(2)
\oplus \mathfrak{sl}(2), \mathfrak{sl}(2)\oplus \mathfrak{sl}(3),
 \mathfrak{osp}(1,2),
  \mathfrak{sl}(2)\oplus \mathfrak{osp}(1,2),$$
$$ \mathfrak{sl}(3)\oplus \mathfrak{osp}(1,2),
\mathfrak{osp}(3,2), \mathfrak{osp}(4,2), \mathfrak{sl}(2)\oplus
\mathfrak{osp}(3,2)\}.$$
  \item[(8)] If $\G$ is of type $F(4)^{\hat{}}$ then
$$\G_{\P}^0\in \{\mathfrak{sl}(2), \mathfrak{sl}(2)
\oplus \mathfrak{sl}(2), \mathfrak{sl}(2)\oplus \mathfrak{sl}(2)
\oplus \mathfrak{sl}(2),  \mathfrak{so}(5), \mathfrak{sl}(2)\oplus \mathfrak{so}(5),$$
 $$\mathfrak{sl}(3),  \mathfrak{sl}(2)\oplus \mathfrak{sl}(3), \mathfrak{sl}(2)\oplus
\mathfrak{sl}(4), \mathfrak{sl}(4), \mathfrak{osp}(4,2),
\mathfrak{sl}(2)\oplus \mathfrak{osp}(4,2)\}.$$
\end{itemize}
 } 

\medskip

The next question we address in the paper is the integrability of
irreducible modules in the category $K(\G)$.  In the case when $\G$ is affine Lie algebra the
classification of all irreducible integrable modules with
finite-dimensional weight spaces was obtained by Chari in
\cite{C}. If the central element acts nontrivially then such
module is highest weight (for some choice of positive roots),
otherwise it is a {\em loop} module
 introduced by Chari and Pressley in \cite{CP}.
 
Analogs of loop modules for affine Lie superalgebras were
constructed by Dimitrov and Penkov in \cite{DP}.
 For $\G$ different
from $A(m,n)^{\hat{}}$ and $C(n)^{\hat{}}$ the irreducible
integrable loop modules were classified Rao and Zhao in
\cite{RZ}. Moreover only trivial integrable modules exist if $\G$ is not of type $A(0,n)^{\hat{}}$,
$B(0,n)^{\hat{}}$, $n\geq 1$, $C(n)^{\hat{}}$, $n\geq 3$ if the central element acts
non-trivially. And if $\G$ has one of these types then irreducible integrable module is highest weight
with respect to some choice of a Borel subalgebra. Integrable highest weight 
modules were classified by Kac and Wakimoto \cite{KW1},
\cite{KW2}. 

It was shown in \cite{KW1} that
the usual definition of integrability for the most of affine Lie
superalgebras allows only trivial highest
 weight modules. Therefore, following \cite{KW1} we introduce  a concept of {\it weak
integrability} for $\G$-modules. Weakly integrable irreducible
highest weight modules of non-twisted affine Lie superalgebras
were classified by Kac and
 Wakimoto in \cite{KW2}. We show that they are the only weakly integrable irreducible
 modules in the category $K(\G)$.
Namely we prove
 the following

\medskip

 \noindent {\bf Theorem C}.
{\it Let $V$ be a weak integrable irreducible module in $K(\G)$. Then $V$ is a highest weight module.}

 \medskip

Theorem C together with the classification of weak integrable
highest weight modules in \cite{KW2} gives a classification of
weak integrable irreducible modules in $K(\G)$.

\section{ Preliminaries}
Let $\U(\G)$ be the universal enveloping algebra of $\G$.
For a subset $S\subset \pi$ denote by $\D(S)$ the additive
closure of $S$ in $\D$, $\D(S)=\Z S\cap \D$.
 Then $\D_+=\D_+(\pi)$
(respectively $\D_-=\D_-(\pi)$) is the set of all positive (resp.,
negative) roots in $\D$ with respect to $\pi$. Set
$\D^{\re}_{\pm}=\D^{\re}\cap \D_\pm$, $\D^{\im}_\pm=\D^\im\cap
\D_\pm =\{k\d|\ k\in \Z_{\pm}\setminus \{0\} \}$, where
 $\d$ is the indivisible positive imaginary root. We will always assume that $\pi$ contains 
only \emph{one} simple odd root of 
$\g$.

For $\a \in \D$ denote by $\G_{\a}$ the root subspace
corresponding to  $\a$.

Let $G = {\C}c \oplus
\displaystyle{\sum_{k \in {\Z} \setminus \{0\} }}
{\G}_{k\d} \subset {\G}$, $G_{\pm} = \displaystyle{\sum_{\pm k>0} {\G}_{k\d}}$.
 For  $a \in
\C$ denote by $M_+(a)$ (resp. $M_-(a)$) the Verma module for $G$ generated by a highest
(lowest) vector $v$ such that $G_+v=0$ (resp. $G_-v=0$) and $cv=av$.
It is known (\cite{F},Proposition 4.5, \cite{FT},Theorem 2.3) that any finitely generated
$\Z$-graded $G$-module $V$ with $\dim V_i<\infty$ for at least one $i\in \Z$ and with
$cV=aV$, $a\neq 0$, is completely reducible. Moreover, all irreducible components
of $V$ are isomorphic $M_+(a)$ or $M_-(a)$ up to a shift of the gradation.

For a $\G$-module $V$ in $K(\G)$ denote $V^{\pm}=\{v\in V|G_{\pm n\delta}v=0
\mbox{\rm\ for all} \, n\geq 1\}$. Following \cite{FT} we define  maximal and
 minimal elements in $V$. A nonzero element $v\in V_{\l}$, $\l\in \H^*$,
is {\it maximal} (resp. {\it minimal}) if $v\in V^+$ (resp. $v\in V^-$)
and $V_{\l+k\delta}\cap V^+=\{0\}$ for all $k>0$ (resp. $V_{\l+k\delta}\cap V^-=\{0\}$ for all $k<0$).

\begin{prop}\label{pr1} Let $V$ be a module in $K(\G)$,
$\l\in P(V)$.
\begin{itemize}
\item [(i)]
The
$G$-module
$V'=\sum_{k\in \Z}V_{\l+ k\d}$
contains a maximal or a minimal element.
\item[(ii)]
Let $v$ be a nonzero  element in $V^+$. Then for any $\a\in{\D}$ there exists $m_0\in\ \Z$ such that
 $\G_{\a+ m\d} v=0$ for all $m\geq m_0$.
\item[(iii)] If $V$ is irreducible then $P(V)$ is a proper subset of a $Q$-coset.
\end{itemize}
\end{prop}

\begin{proof} Statement (i) is an analog of Proposition 3.1 in \cite{FT}.
Consider the
$G$-module
$V''$
of $V'$ generated by
$V_\l$.
Then the module
$V''$
is completely reducible  and all its irreducible
submodules
are of type
$M_{\pm}(\l(c))$ up to the shift of gradation.
Suppose that
$V^+\cap V'\not=0$.
Since
$(V'/V'')_\l=0$, then
$V^+\cap V_{\l+k\d}\subset V''$
for all
$k\geq 0$ and there exists only a finite number of integers
$k\geq 0$ for which $V^+\cap V_{\l+k\d}\not=0$. Hence $V'$ contains a maximal element.
 If $V^+\cap V'=0$
then clearly $V^-\cap V'\not=0$ and the same arguments as above show
that $V'$ contains a  minimal element.

Statement (ii) generalizes Lemma 3.2 in \cite{FT}.
 Assume that $\a \in \D^\re$. Let
$v\in V_\l \cap V^+$. For each $\f\in \D^\re$ fix a nonzero element
$X_{\f}\in \G_{\f}$. Suppose that $w_m=X_{\a+m\d}v\not=0$ for all
$m\in \Z_+$.  Since
$[\G_{m\d},\G_\a]=\G_{\a+m\d}$, then
 a $G$-submodule $N$ generated by $w_0$ contains all $w_m$, $m\in \Z_+$. The module $N$
is completely reducible with irreducible submodules isomorphic to $M_{\pm}(\l(c))$. Since
$\dim V_\l<\infty$,
 there exists $m_0 \in \Z_+$ such that
the $G_+$-submodule $\U(G_+)w_{m_0}$ is $\U(G_+)$-free. On the other hand,
$$[\G_{2\d}, \G_{\a+m_0 \d}]=
[\G_{\d}, [\G_{\d}, \G_{\a+m_0 \d}]]=\G_{\a+(m_0+2)\d}
$$
and hence $x w_{m_0}$ is proportional to the element $y^2w_{m_0}$
for all $x\in \G_{2\d}$,
$y\in \G_{\d}$ since the real root spaces are $1$-dimensional.
 Since $v\in V^+$ then the obtained contradiction shows that $\G_{\a+m\d} v=0$
for all $m$ sufficiently large. This completes the proof of (ii).

Statement (iii) is an analog of Proposition 3.4 in \cite{FT}. Suppose $V$ is irreducible.
Without loss of generality we may assume by (i) that $V$ contains a maximal element
$v\in V_{\l}\setminus \{0\}$. Let $\mu\in P(V)$.
Since $V$ is irreducible, there exists
$u\in \U(\G)$ such that $0\not= uv\in V_\mu$.
 Applying (ii) we conclude that $\G_{n\d}uv=0$
for sufficiently large $n$. From here we deduce  that a subspace $V'=\sum_{k\in \Z_+}V_{\mu+k\d}$
has a nonzero intersection with $V^+$. The $G$-submodule $\U(G)uv$
is completely reducible with
irreducible components of type
$M_{\pm}(\l(c))$, moreover there exists an irreducible component
$N$ of type $M_+(\l(c))$ such that $N\cap V_\mu\not=0$. As in the proof of (i) we can show
 the existence of a maximal element in $V'$. Now suppose that $n>0$ and $\lambda +n\d\in P(V)$.
 Applying the same argument as above to the weight
$\mu=\l+n\d$, we conclude that
$$V^+\cap(\sum_{k\in \Z_+}V_{\l+(k+n)\d})\not=0.$$
But this contradicts our assumption on the maximality of
$v$ which shows that $\l+n\d\notin P(V)$ for all $n>0$. Hence $P(V)$ is a proper subset of a $Q$-coset.
\end{proof}

\begin{remark}
Proposition~\ref{pr1} also holds for $\mathfrak{psl}(n,n)$. Indeed all arguments valid for $n>2$. If $n=1$ then all root spaces are $2$-dimensional. In this case the proof of item (ii) can be modified as follows. Consider the element 
$w_{m_0}$ from the proof of  Proposition~\ref{pr1},(ii). For each $m>0$ choose a nonzero $x_m\in \G_{m\delta}$. Then for any positive $s$ the vectors
$$\{x_{i_1}\ldots x_{i_r}w_{m_0}|i_1+\ldots +i_r=s, i_j>0, j=1, \ldots, r,\}$$
are linearly independent. On the other hand
$$[\G_{i_r\d}, [\ldots [\G_{i_1\d}, \G_{\a+m_0 \d}]\ldots] ]=\G_{\a+(m_0+s)\d},$$
and hence the space spanned by all $x_{i_1}\ldots x_{i_r}w_{m_0}$, $i_1+\ldots +i_r=s$, is at most $2$-dimensional, which is a contradiction. 
 
\end{remark}

Without  loss of generality we will always assume that modules of $K(\G)$ contain
a maximal element.

\section{Cuspidal Levi subsuperalgebras}
Let $\P$ be a parabolic subset of $\D$,
 $\G=\G_{\P}^-\oplus (\G_{\P}^0+\H)\oplus \G_{\P}^+$ where
$\G_{\P}^{\pm}=\sum_{\a\in \P}\G_{\pm \a}$ and  $\G_{\P}^0$
is generated by the subspaces
 $\G_{\a}$ with $\a\in {\P}^0$.
 The subsuperalgebra $\G_{\P}^0$ is called a {\em Levi subsuperalgebra} of $\G$.
Let $N$ be an irreducible weight module over a parabolic subalgebra
 $\G_{\P}=(\G_{\P}^0+\H)\oplus \G_{\P}^+$, with a
trivial action of
$\G_{\P}^+$ and let
$$M_{\P}(N)={\rm ind}(\G_{\P}, \G; N)$$
 be the induced $\G$-module. Module $M_{\P}(N)$ has a unique irreducible quotient
$L_{\P}(N)$. If $\G_{\P}^0\neq \G$ then $L_{\P}(N)$ is said to be
{\em parabolically induced}.
An irreducible weight $\G$-module is called {\em cuspidal} if is not parabolically induced.
 A Levi subsuperalgebra of $\G$ is {\em cuspidal} if it admits a weight cuspidal module with
finite-dimensional weight spaces. All cuspidal Levi subalgebras of reductive Lie algebras were
classified by Fernando \cite{Fe}. They are the subalgebras with simple components of type $A$
and $C$. Cuspidal Levi subsuperalgebras of finite-dimensional Lie superalgebras  were
 described by Dimitrov, Mathieu and Penkov in \cite{DMP} and \cite{DMP1}.
  Here we present their  list.
Besides the  cuspidal Levi subalgebras
 of type $A$ and $C$ it includes the superalgebras of type $\mathfrak{osp}(1,2)$,  $\mathfrak{osp}(1,2)\oplus
 \mathfrak{sl}(2)$,
 $\mathfrak{osp}(n,2m)$
with $2<n\leq 6$ and  $D(2,1;a)$. Note that notion a
of cuspidal
Levi subsuperalgebra
in \cite{DMP} is more general  than ours since it refers to generalized weight modules. We will classify all
cuspidal  subsuperalgebras of affine Lie superalgebras in Section 6.

\section{Parabolically induced modules}

Let $\P \subset \D$ be a parabolic subset.
A nonzero element $v$ of a $\G$-module $V$ is called $\P$-{\it primitive} if $\G_{\P}^+v=0$.
  The module $M_{\P}(N)$ has a $\G_{\P}^0\oplus \G_{\P}^+$-submodule isomorphic to $N$ which consists of
$\P$-primitive elements.

Note that the category $K(\G)$ is closed under taking the submodules and the quotients.
If $N$ has a finite-dimensional weight space and $c$ acts on $N$ as a nonzero scalar
then modules $M_{\P}(N)$ and $L_{\P}(N)$ are objects of $K(\G)$.
Let $Q_{\P}$ be the free abelian group generated by $\P^0$. The universality of the module
$M_{\P}(N)$ is clear from the following standard statement (cf. \cite{FT},Proposition 2.2).

\begin{prop}\label{pr2} Let $\P^0\neq \pi$, $V$ be an irreducible weight $\G$-module
with a $\P$-primitive element of weight $\l$,
 $N=\sum_{\nu \in Q_{\P}}V_{\l+\nu}$. Then  $N$ is an
irreducible $\G_{\P}^0\oplus \G_{\P}^+$-module and $V$ is isomorphic to $L_{\P}(N)$.
\end{prop}
\bigskip

\begin{remark}
The module $N$ in Proposition~\ref{pr2} is an irreducible weight module over
a finite-dimensional Lie superalgebra $\G_{\P}^0$. By \cite{DMP},Theorem 6.1, any such module is
parabolically induced from a cuspidal module over a subsuperalgebra of $\G_{\P}^0$. Hence the
classification of all irreducible weight $\G$-modules with a $\P$-primitive element
 is reduced to the classification of cuspidal modules over finite-dimensional Lie superalgebras.
\end{remark}

If a parabolic subset $\P$ is standard (with respect to $\pi$),
i.e. $\P^+\subset \D_+$,  $\P^0=\D(S)$ for some $S\subset \pi$, then we write
$L_{\P}(N)=L^S(N)$.  
If $\cB\subset\G$ is a subsuperalgebra then
we will also consider irreducible modules $L^T_{\cB}(K)$ obtained
in the same manner, where $T$ is a subset of a basis of the root
system of $\cB$ and $K$ is an irreducible
 $\cB(T)$-module.

A $\P$-primitive element with $\P=\P_S$ will be called
$S$-primitive. If $V$ is generated by a $\emptyset$-primitive element $v\in V_{\l}$
then $V$ is a highest weight module with highest weight $\l$.

The situation is especially pleasant when $S\subset \pi_0$. In
this case $\G(S)$ is a finite-dimensional semisimple Lie algebra
and we have the following stronger version of
Proposition~\ref{pr2}.

\begin{prop}\label{pr3}
Let $S\subset \pi_0$, $V$ be an irreducible weight $\G$-module with a $S$-primitive element.
Then there exists a basis $\pi'$ of $\D$, a subset $S'\subset \pi'_0$ and an irreducible cuspidal
$\G(S')+\H$-module $N$ such that $V\simeq L^{S'}(N)$. Moreover, if $V$ has finite-dimensional
weight spaces then all simple components of $\G(S')$ are of type $A$ and $C$.
\end{prop}

In particular, Proposition~\ref{pr3} together with \cite{M}
describes all irreducible modules with $S$-primitive elements in
the category $K(\G)$. Moreover, in the case when $\G$ is an affine
Lie algebra the modules of type $L^S(N)$ with cuspidal $N$ exhaust
all irreducible modules in $K(\G)$. Namely, the following result
was proved in \cite{FT}.

\begin{theorem}\label{thm1}
Let $\G$ be an affine Lie algebra, $V$ a $\G$-module in $K(\G)$.
\begin{itemize}
\item[(i)] (\cite{FT},Proposition 4.3) There exists a basis $\pi$ of $\D$ and $\alpha\in \pi$
such that $V$ contains a $S$-primitive element, $S=\pi\setminus \{\alpha\}$.
\item[(ii)] (\cite{FT},Theorem 4.1) If $V$ is irreducible then there exists a basis $\pi$
of $\D$, a proper subset $S\subset \pi$ and an irreducible weight cuspidal $\G(S)+\H$-module
 $N$ such that $V\simeq L^S(N)$.
\end{itemize}
\end{theorem}

\section{Existence of $\P$-primitive elements in objects of $K(\G)$}

Let $\D=\D_0\cup \D_1$ where $\D_0$ (resp. $\D_1$) consists of all even (resp. odd)
roots, $\D^{\pm}_0=\D_0\cap \D_{\pm}$, $\D^{\pm}_1=\D_1\cap \D_{\pm}$.
 Recall that $\G=\G_0\oplus \G_1$  and
$\g=\g_0\oplus \g_1$.
 Here $\g_0$ is a reductive Lie algebra with at most three simple
 components. In this section we will denote by $\cB_0$ the semisimple part of
 $\g_0$ and by $\cB$ its affinization.

  Let $\dot{\D}$ be the root system of $\g$,
$\dot{\D}=\dot{\D}_0\cup\dot{\D}_1$, $\dot{\D}_0^{\pm}=
\dot{\D}_0\cap \D_{\pm}$, $\dot{\D}_1^{\pm}= \dot{\D}_1\cap
\D_{\pm}$, $\dot{\D}_+=\dot{\D}_0^+ \cup \dot{\D}_1^+$,
 $\dot{\pi}$ a basis of $\dot{\D}$. Then
$\D^{\re}=\dot{\D}+\Z\d$, $\D^{\re}_+=\{\dot{\D}_+ +\Z_{\geq
0}\d\}\cup \{-\dot{\D}_+ +\Z_{> 0}\d\}$,
$\D_0=(\dot{\D}_0+\Z\d)\cup \D^{\im}$. Let $W_0$ be the Weyl group
of $\G_0$.

In this section we establish the following key result for the proof of Theorem A.

\begin{prop}\label{pr4}
Let $V$ be a nonzero $\G$-module in $K(\G)$. Then $V$ contains a $\P$-primitive element
for some proper parabolic subset $\P$.
\end{prop}

We start with the following 
\begin{lemma}\label{lemma6} Let $V$ be a nonzero module in $K(\G)$ and
$V_{\cB}$ the restriction of $V$ on $\cB$. Let also
$\cB_0=\cB_0^1\oplus \cB_0^2$ and  denote by $\cB^1$ and $\cB^2$
the corresponding affine Lie algebras. Then there exists a basis
$\pi_{\cB}=\pi_{\cB}^1\cup\pi_{\cB}^2$ of the root system of $\cB$
and $\a_i\in \pi_{\cB}^i$ such that $V_{\cB}$ contains a
$S$-primitive element $v\in V_{\cB}$, $S=\pi_{\cB}\setminus
\{\a_1, \a_2\}$, where $\pi_{\cB}^i$ is  a basis of the root
system of $\cB_i$,   $i=1,2$.
\end{lemma}

\begin{proof}   Let
$W_{i}\subset W_0$ be the Weyl group. Consider $V_{\cB}$ as a
$\cB^1$-module and apply Theorem~\ref{thm1}(i). Then there
exists $w\in W_{1}$, a basis $\tilde{\pi}_{\cB}^1=w{\pi}_{\cB}^1$
and a root $\a_1\in \tilde{\pi}_{\cB}^1$ such that $V$ contains a
$S_1$-primitive element
 $v'$ for $\cB_0^1$, $S_1=\tilde{\pi}_{\cB}^1\setminus \{\a_1\}$.
 Consider a basis $w\pi$ of $\D$, note that
$w\pi_{\cB}^2=\pi_{\cB}^2$. Denote by $V'$ a subspace of $V^+$
consisting of $S_1$-primitive elements  for $\cB_0^1$, $v'\in V'$.
Clearly, $V'$ is a $\cB_0^2$-module. Applying again
Theorem~\ref{thm1}(i)  to $V'$ we find  $w'\in W_{2}$, a basis
$\tilde{\pi}_{\cB}^2=w'{\pi}_{\cB}^2$ and a root $\a_2\in
\tilde{\pi}_{\cB}^2$ such that $V'$ contains a $S_2$-primitive
element
 $v$ for $\cB_0^2$, $S_2=\tilde{\pi}_{\cB}^2\setminus \{\a_2\}$.
 Then $\tilde{\pi}_{\cB}^1\cup\tilde{\pi}_{\cB}^2$ is a required
 basis of the root system of $\cB$ and $v$ is $S=S_1\cup S_2$-primitive element.
\end{proof}

Let $V$ be a $\G$-module in $K(\G)$. For $v\in V\setminus \{0\}$ denote by $\Omega_v$ the set of
all roots $\alpha\in \D_+$ such that $\G_{\alpha}v=0$. In particular, $\D_+^{\im}\subset \Omega_v$
if $v\in V^+$. Moreover, in this case $\alpha+N\d\in \Omega_v$ for any $\alpha\in \dot{\D}$ and
any $N$ sufficiently large by Proposition~\ref{pr1}(ii).

Given $\alpha\in \D^{\re}$ denote
$$\Sigma_{\alpha}=\{\alpha+n\d, n\in \Z_{\geq 0}\}\cup \{-\alpha+m\d, m\in \Z_{>0}\}.$$
If $\alpha\in \D_0^{+}$ then $\Sigma_{\alpha}$ consists of positive real roots of the affine
$\mathfrak{sl}(2)$-subalgebra of $\G$ corresponding to the root $\alpha$. Note that
$\Sigma_{\alpha}=\Sigma_{-\alpha+\d}$. If $V$ is a $\G$-module then an element $v\in V^+\setminus \{0\}$
 is said to be $(\a,\pi)$-{\it admissible} if $\Sigma_{\a}\setminus \{\a\}\subset \Omega_v$.
If $S\subset \dot{\pi}$ then it follows immediately from the definition that $v\in V^+$
is $S$-primitive \iff it is $(\a,\pi)$-admissible for all $\a\in \D(S)\cap \D_+$ and
 $(-\beta+\d, \pi)$-admissible for all $\beta\in \dot{\D}_+\setminus \D(S)$. We will say
that $v\in V^+$ is $\pi$-admissible if it is  $(\b,
\pi)$-admissible or $(-\b+\d, \pi)$-admissible for all $\b\in
\dot{\D}_0^+$.

Recall that an odd simple root $\ga\in \dot{\pi}$ is isotropic if $(\ga,\ga)=0$
which is if and only if $[\G_{\ga}, \G_{\ga}]=0$.

\begin{lemma}\label{le1}
\begin{itemize}

\item[(i)] Let $V$ be a weight $\G$-module that contains
 a $\pi$-admissible element. Then there
exists a $\pi$-admissible $v\in V^+$ which is $(\ve,
\pi)$-admissible for all $\ve\in \dot{\D}^+_1$. Moreover, if $v$
is $(\b, \pi)$-admissible for all $\b\in \dot{\D}^+_0$ and
$[\G_{\ve}, \G_{\ve'}]=0$ for all $\ve, \ve'\in \dot{\D}^+_1$ then
$V$ contains a $\dot{\pi}\setminus \{\gamma\}$-primitive element, where
$\ga$ is the odd simple root in $\dot{\pi}$.

\item[(ii)] Let $\ga\in \dot{\pi}$ be an isotropic root. If
$[\g_{\ve}, \g_{\ve'}]=0$ for all $\ve, \ve'\in \dot{\D}^+_1$ then
any weight $\g$-module $V$ contains a $S$-primitive element,
$S=\dot{\pi}\setminus \{\ga\}$.

\item[(iii)]  Let $V$ be a weight $\G$-module and $v$ a nonzero element of $V$ such that
$\D_0^+\subset \Omega_v$. Then $V$ contains
a $\P$-primitive element, where $\P=\P_{\emptyset}=\D_+$, thus $V$ is a highest weight module.

\end{itemize}
\end{lemma}

\begin{proof}
Let $V'\subset V$ be a subspace of $\pi$-admissible elements in
$V$ and $v\in V'\setminus \{0\}$. By Proposition~\ref{pr1}(ii)
there exists $N>0$ such that $\ve+n\d\in \Omega_v$ for all $n\geq
N$ and all $\ve\in \dot{\D}^+_1$. Let $N_0$ be the least
nonnegative integer such that $\ve+N_0\d\in \Omega_v$ for all
$\ve\in \dot{\D}^+_1$ and there exists $\ve'\in \dot{\D}^+_1$ such
that $\ve'+(N_0-1)\d\notin \Omega_v$. Suppose $N_0>1$ and choose
 $\ve_0\in \dot{\D}^+_1$  of maximal height with $\ve_0+(N_0-1)\d\notin \Omega_v$. Consider
$w=X_{\ve_0+(N_0-1)\d}v\neq 0$ which obviously belongs to $V'$. Moreover, $\ve_0+(N_0-1)\d\in \Omega_w$
and hence $\ve_0+k\d\in \Omega_w$ for all $k\geq N_0-1$. Since $\dot{\D}_1$ is finite we can find
 a nonzero $v'\in V'$ such that $\ve+k\d\in \Omega_{v'}$ for all $\ve\in \dot{\D}^+_1$ and $k>0$.
Considering now the roots $-\ve+n\d$,  $\ve\in \dot{\D}^+_1$, $n>0$, we choose $N_0$ as above
 and $\ve_0$ of minimal height such that  $\ve_0+(N_0-1)\d\notin \Omega_{v'}$. Then the same argument
as above show the existence of $w\in V'\setminus \{0\}$ such that $\a+k\d\in \Omega_w$ for all
$\a\in \dot{\D}^+_1$ and $k>0$. Hence $w$ is $(\ve, \pi)$-admissible for all $\ve\in \dot{\D}^+_1$.

Suppose now that  $v$ is $(\b, \pi)$-admissible for all $\b\in \dot{\D}^+_0$ and
$[\G_{\ve}, \G_{\ve'}]v=0$ for all $\ve, \ve'\in \dot{\D}^+_1$. Then in particular
$v$ is $\dot{\pi}$-primitive. Moreover, in this case using the same
arguments as above we find a nonzero
$w\in V'$ which is  $(\ve, \pi)$-admissible for all $\ve\in \dot{\D}^+_1$ and
$\dot{\D}^+_1\subset \Omega_w$.  Therefore $\G_{\b}w=0$ for all $\b\in \D_+\setminus
\D(\dot{\pi}\setminus \{\ga\})$ which proves (i).
The proof of (ii) and (iii) is similar.
\end{proof}

\begin{remark}
 The idea of the proof of Lemma~\ref{le1}  will be used in
 various cases in this section. We will leave the details out since the
 exposition is similar to the one above.
\end{remark}

We fix  a nonzero $\G$-module $V$ in $K(\G)$.

\subsection{Case  $A(m,n)^{\hat{}}$}

Let $\pi_{01}=\{\a_1, \ldots, \a_n\}$, $\pi_{02}=\{\b_1, \ldots,
\b_m\}$, $\pi_0=\pi_{01}\cup \pi_{02}$, $\pi_1=\{\ga, \a_0\}$
according to the following diagram

\bigskip

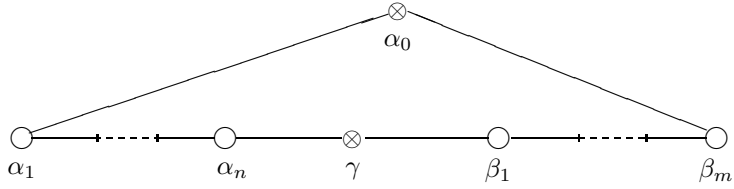
\unitlength=1mm
\special{em:linewidth 0.4pt}
\linethickness{0.4pt}
\begin{picture}(19.75, 30.41)
\put(  9.00,11.67){\line(3,1){47.67}}
\put( 9.30,11.00){\line(1,0){8.70}}
\put( 26.00,11.00){\line(1,0){7.60}}
\put( 36.50,11.00){\line(1,0){14.00}}
\put( 70.00,11.00){\line(-1,0){16.30}}
\put( 73.00,11.00){\line(1,0){9.00}}
\put( 99.00,12.00){\line(-5,2){39.67}}
\put( 91.00,11.00){\line(1,0){8.00}}
\put(  8.00,11.00){\circle{2.83}}
\put( 35.00,11.00){\circle{2.83}}
\put( 56.55,27.00){$\otimes$}

\put( 71.33,11.00){\circle{2.83}}
\put(100.33,11.00){\circle{2.83}}
\put( 82.00,11.00){\dashbox{0.8}(9.00,0.00)[cc]{}}
\put( 18.00,11.00){\dashbox{0.8}(8.00,0.00)[cc]{}}
\put(  8.00, 6.67){\makebox(0,0)[cc]{$\alpha_1$}}
\put( 52.00, 6.67){\makebox(0,0)[cc]{$\gamma$}}
\put( 36.00, 6.67){\makebox(0,0)[cc]{$\alpha_n$}}
\put( 58.00,23.67){\makebox(0,0)[cc]{$\alpha_0$}}
\put( 71.33, 6.67){\makebox(0,0)[cc]{$\beta_1$}}
\put(100.33, 6.67){\makebox(0,0)[cc]{$\beta_m$}}
\put( 50.55,10.00){$\otimes$}

\end{picture}

Let $\g_{01}$ and $\g_{02}$ be the corresponding simple components of $\g_0$ with the root systems
$\dot{\D}_{01}$ and $\dot{\D}_{02}$ respectively. Set $\theta=\a_1+\ldots + \a_n$.

\begin{lemma}\label{le2}
\begin{itemize}
\item[(i)] For any $\b\in \pi_{01}$ there exists a basis $\pi'$ of $\D$ such that
$\pi'_0=\pi'_{01}\cup \pi'_{02}$, $\pi'_{02}=\pi_{02}$,
$\pi'_{01}=(\pi_{01}\setminus \{\b\})\cup \{-\theta+\d\}$.
\item[(ii)] Let $\b\in \pi_{01}$. If $v$ is $(\a, \pi)$-admissible for all $\a\in
\D(\pi_{01}\setminus \{\b\})\cap \D^+_0$ and $(-\ve+\d, \pi)$-admissible for all
$\ve\in \dot{\D}^+_{01}\setminus \D(\pi_{01}\setminus \{\b\})$ then there exists a basis
$\pi'$ of $\D$ such that $\pi'_0=\pi'_{01}\cup \pi'_{02}$, $\pi'_{02}=\pi_{02}$
and $v$ is $(\a, \pi')$-admissible for
$\a\in \D(\pi'_{01})\cap \D_0^+(\pi')$.
\end{itemize}
\end{lemma}

\begin{proof}
For any $i=1, \ldots, n$ there exists a new basis $\pi'$ of $\D$ such that $\pi'_0=\pi'_{01}\cup \pi'_{02}$,
 $\pi'_{01}=\{\a'_1, \ldots, \a'_n\}$ where $\a'_k=\a_{i-k}$, $k=1, \ldots, i-1$, $\a'_i=-\theta+\d=
\ga+\b_1+\ldots +\b_m+\a_0$, $\a'_{i+k}=\a_{n-k+1}$, $k=1, \ldots, n-i$, $\pi'_{02}=\{\b_1, \ldots, \b_m\}$,
$\pi'_1=\{-\a_{i+1}-\ldots -\a_n-\ga-\b_1-\ldots -\b_m, \ga+\a_n+ \ldots+ \a_i\}$ which implies (i).
Given $\b\in\pi'_{01}$ choose a new basis $\pi'$ such that $\pi'_{01}=(\pi_{01}\setminus \{\b\})\cup
\{-\theta+\d\}$, $\pi'_{02}=\pi_{02}$ which exists by (i). Since $\D(\pi_{01}\setminus \{\b\})\subset
\D(\pi'_{01})$ and $-\ve+\d\in \D(\pi'_{01})$ for all $\ve\in \dot{\D}^+_{01}\setminus
\D(\pi_{01}\setminus \{\b\})$ we conclude that $\pi'$ is a required basis which proves (ii).
\end{proof}

Denote by $\G_{01}$ and $\G_{02}$ the affine Lie algebras with underlying finite-dimensional
algebras $\g_{01}$ and $\g_{02}$, and with the  root systems $\D_{01}$ and $\D_{02}$ respectively.
 Then $\tilde{\pi}_{01}=\pi_{01}\cup \{-\theta +\d\}$ is a basis of $\D_{01}$. Let $W_{0i}\subset W_0$
be the Weyl group of $\D_{0i}$, $i=1,2$. As in the proof of
Lemma~\ref{lemma6} consider $V$ as a $\G_{01}$-module and apply
Theorem~\ref{thm1}(i). Then there exists $w\in W_{01}$, a basis
$\tilde{\pi}'_{01}=w\tilde{\pi}_{01}$ of $\D_{01}$ and a root
$\a\in \tilde{\pi}'_{01}$ such that $V$ contains a $S$-primitive
element
 $v'\in V$, $S=\tilde{\pi}'_{01}\setminus \{\a\}$. Consider a basis $w\pi$ of $\D$, note that
$w\pi_{02}=\pi_{02}$. Applying Lemma~\ref{le2}(ii) to $\a$ we conclude that there exists a basis
$\pi'$ of $\D$ such that $\pi'_0=\pi'_{01}\cup \pi_{02}$ and $v'$ is $(\b, \pi')$-admissible
for all $\b\in \D(\pi'_{01})\cap w\D_0^+$. Denote by $V'$ a subspace of $V^+$ consisting
of elements with such property ($v'\in V'$). Clearly, $V'$ is a $\G_{02}$-module.
Applying subsequently Theorem~\ref{thm1}(i) and  Lemma~\ref{le2}(ii) to $V'$ we conclude that
there exists a basis $\pi''$ of $\D$ and a nonzero $\tilde{v}\in V'$ which is
$(\a, \pi'')$-admissible for all $\a\in \D_+(\pi_0'')$. Let $\tilde{V}$ be a subspace of $V'$
consisting of  elements with such property and apply Lemma~\ref{le1}(i) to $\tilde{V}$. Since
 $[\G_{\ve}, \G_{\ve'}]=0$ for all $\ve, \ve'\in \dot{\D}(\pi'')^+_1$, there exists a nonzero
 $\bar{v}\in \tilde{V}$ which is $\pi''_0$-primitive. If $\gamma\in \pi''$ is an odd root then
in particular, $\bar{v}$ is $S$-primitive, $S=\pi''\setminus
\{\ga\}$ which completes the proof of Proposition~\ref{pr4} for
$A(m,n)^{\hat{}}$.

\medskip

Consider the restriction $V_{\cB}$ of $V$ on $\cB=\sum_{i\in I}
\cB_i$. Then by Lemma~\ref{lemma6}  there exists a basis
$\pi_{\cB}=\cup_{i\in I} \pi_{\cB}^i$ of the root system of $\cB$
and $\a_i\in \pi_{\cB}^i$ for each $i\in I$ such that $V_{\cB}$
contains a $S_{\cB}$-primitive element $v\in V_{\cB}$,
$S_{\cB}=\pi_{\cB}\setminus \{\a_i, i\in I\}$. In particular, $v$ is
$\pi'$-admissible where $\pi'$ is a basis of $\D$ containing
$\pi_{\cB}$. In all subsequent cases we  fix
 such   element $v\in V_{\cB}$ and assume without loss of generality
  that $\pi'=\pi$.

\subsection{Case $C(n)^{\hat{}}$,  $n\geq 3$}

Let $\G$ be of type $C(n)^{\hat{}}$, $\pi=\{\a_0, \a_1, \ldots, \a_n\}$
a basis of $\D$ ordered as follows

\bigskip

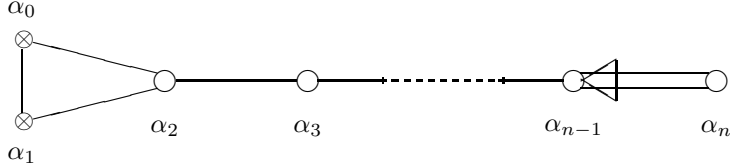
\unitlength=1.00mm
\special{em:linewidth 0.4pt}
\linethickness{0.4pt}
\begin{picture}(101.74,23.80)
\put(28.50,11.00){\line(1,0){16.00}}
\put(27.00,11.00){\circle{2.83}}
\put(46.00,11.00){\circle{2.83}}
\put(81.33,11.00){\circle{2.83}}
\put(100.33,11.00){\circle{2.83}}
\put(47.33,11.00){\line(1,0){9.00}}
\put(80.00,11.00){\line(-1,0){8.67}}
\put(56.00,11.00){\dashbox{0.8}(16.00,0.00)[cc]{}}
\put(8.00,20.50){\makebox(0,0)[cc]{$\alpha_0$}}
\put(8.00,1.00){\makebox(0,0)[cc]{$\alpha_1$}}
\put(27.00,4.67){\makebox(0,0)[cc]{$\alpha_2$}}
\put(81.33,4.67){\makebox(0,0)[cc]{$\alpha_{n-1}$}}
\put(100.33,4.67){\makebox(0,0)[cc]{$\alpha_n$}}
\put(46.00,4.67){\makebox(0,0)[cc]{$\alpha_3$}}
\put(82.30,10.00){\line(1,0){17.00}}
\put(82.30,12.00){\line(1,0){17.00}}
\put(26.00,12.00){\line(-4,1){16.5}}
\put(26.00,10.00){\line(-4,-1){16.5}}
\put(82.33,11.00){\line(5,3){4.67}}
\put(87.00,13.67){\line(0,-1){5.33}}
\put(8.00,15.10){\line(0,-1){8.33}}
\put(87.00,8.33){\line(-5,3){4.67}}
\put(6.87,4.74){$\otimes$}
\put(6.87,15.58){$\otimes$}
\end{picture}

\bigskip

Here $\dot{\pi}=\{\a_1, \ldots, \a_n\}$ is a basis of $\dot{\D}$
and $\pi_0=\{\a_2, \ldots, \a_n\}$ is a basis of $\dot{\D}_0$.
Note that  $\a_0+\a_1\in \D_0\setminus \dot{\D}_0$. Then $\pi_0$
is a basis of  $\cB_0$ which has type $C_{n-1}$, and
$\tilde{\pi}=\pi_0\cup \{-\theta+\d\}$ is a basis of $\cB$ of type
$C_n^{(1)}$, where $\theta$ is the highest positive root in
$\D(\pi_0)$.

By Lemma~\ref{lemma6}, $v$ is a $S_{\cB}$-primitive element, where
$S_{\cB}=\tilde{\pi}\setminus \{\a\}$ for some $\a\in \tilde{\pi}$.
Hence we have the following cases.

\noindent{\it Case 1.} Suppose $\a=-\theta+\d$. Then $v$  is $(\b,
\pi)$-admissible for any $\b\in \D_+(\pi_0)$. Since $[\G_{\ve},
\G_{\ve'}]=0$ for all $\ve, \ve'\in \dot{\D}^+_1$, we can apply
Lemma~\ref{le1}(i)
 to find a $\pi_0$-primitive element in $V$.

\noindent{\it Case 2.} Suppose now that $\a=\a_n$. Hence  $v$  is
$(\b, \pi)$-admissible for all $\b\in \D_+((\pi_0\setminus
\{\a_n\})\cup\{-\theta+\d\})$.
 Choose a new basis $\pi'=\{\a_0',\a_1', \ldots, \a_n'\}$ where $\a_0'=-\a_1-\ldots -\a_{n-1}$,
$\a_1'=\a_1+\ldots +\a_n$, $\a_k'=\a_{n-k+1}$, $k=2, \ldots, n-1$, $\a_n'=\a_0+\a_1$. Then $v$ is
$(\b, \pi')$-admissible for any $\b\in \D_+(\pi_0')$ and we are back to Case 1.

\noindent{\it Case 3.} Let $\a=\a_k$, $2\leq k\leq n-1$. Hence $v$
is
 $(\b, \pi)$-admissible for any $\b\in \D_+(\pi_0\cup \{\a_0+\a_1\})\setminus
(\D_+(\pi_0)\setminus \D_+(\pi_0\setminus \{\a\}))$. Denote a
subspace of such elements by $V'$, hence $V'\neq 0$. Consider the
set $\D^+_{1,\a}$ of positive  roots in $\D_1$ that have $\a$ in
their decomposition in simple roots. Then applying the same
argument as in the proof of Lemma~\ref{le1}(i), we find a nonzero
 $v'\in V'$ such that $\ga\in \Omega_{v'}$ for all
$\ga\in \D^+_{1,\a}$ and $k>0$. Hence $v'$ is $S$-primitive, where
$S=\pi\setminus \{\a\}$. This completes the proof of
Proposition~\ref{pr4} for $C(n)^{\hat{}}$.

\subsection{Case $B(m,n)^{\hat{}}, m>0$}
Let  $\G$ be of type $B(m,n)^{\hat{}}$, $m>0$,  $\pi=\{\a_0,
\ldots, \a_n,$  $\b_1, \ldots, \b_m\}$, $\dot{\pi}_0=\{\a_1,
\ldots, \a_{n-1}, \b_1, \ldots, \b_m\}$, $\dot{\pi}_1=\{\a_n\}$,
$\d=\a_0+2\sum_{i=1}^n +2\sum_{j=1}^m \b_j$.

\medskip

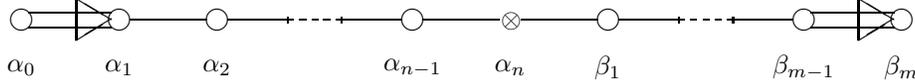
\unitlength=1.00mm
\special{em:linewidth 0.4pt}
\linethickness{0.4pt}
\begin{picture}(120.74,26.80)
\put(1.00,11.00){\circle{2.83}}
\put(14.00,11.00){\circle{2.83}}
\put(2.00,12.00){\line(1,0){11.00}}
\put(2.00,10.00){\line(1,0){11.00}}
\put(27.00,11.00){\circle{2.83}}
\put(53.00,11.00){\circle{2.83}}
\put(64.78,10.00){$\otimes$} \put(79.00,11.00){\circle{2.83}}
\put(105.00,11.00){\circle{2.83}}
\put(118.00,11.00){\circle{2.83}}
\put(15.40,11.00){\line(1,0){10.20}}
\put(54.40,11.00){\line(1,0){10.20}}
\put(67.40,11.00){\line(1,0){10.20}}
\put(28.40,11.00){\line(1,0){8.00}}
\put(106.00,12.00){\line(1,0){11.00}}
\put(106.00,10.00){\line(1,0){11.00}}
\put(43.60,11.00){\line(1,0){8.00}}
\put(88.40,11.00){\dashbox{0.8}(7.20,0.00)[cc]{}}
\put(80.40,11.00){\line(1,0){8.00}}
\put(95.60,11.00){\line(1,0){8.00}}
\put(36.40,11.00){\dashbox{0.8}(7.20,0.00)[cc]{}}
\put(1.00,4.67){\makebox(0,0)[cc]{$\alpha_0$}}
\put(14.00,4.67){\makebox(0,0)[cc]{$\alpha_1$}}
\put(27.00,4.67){\makebox(0,0)[cc]{$\alpha_2$}}
\put(53.00,4.67){\makebox(0,0)[cc]{$\alpha_{n-1}$}}
\put(66.00,4.67){\makebox(0,0)[cc]{$\alpha_n$}}
\put(79.00,4.67){\makebox(0,0)[cc]{$\beta_1$}}
\put(105.00,4.67){\makebox(0,0)[cc]{$\beta_{m-1}$}}
\put(118.00,4.67){\makebox(0,0)[cc]{$\beta_m$}}
\put(117.00,11.00){\line(-5,3){4.67}}
\put(112.30,13.67){\line(0,-1){5.33}}
\put(112.30,8.33){\line(5,3){4.67}}
\put(13.00,11.00){\line(-5,3){4.67}}
\put(8.30,13.67){\line(0,-1){5.33}}
\put(8.30,8.33){\line(5,3){4.67}}
\end{picture}

\medskip

Since  $\cB_0\simeq C_n\oplus B_m$ then  $\cB\simeq \mathfrak{A}_1
+ \mathfrak{A}_2$, where $\mathfrak{A}_1$ is an affine Lie algebra
of type $C_n^{(1)}$ while $\mathfrak{A}_2$ is an affine Lie
algebra of type $B_m^{(1)}$ for $m>2$, of type $C_2^{(1)}$ for
$m=2$ and of type $A_1^{(1)}$ for $m=1$. Denote
$\gamma_1=2\a_n+2\sum_{j=1}^m \b_j$ $\gamma_2=-\b_1-2\sum_{j=2}^m
\b_j+\d$. Then $\{\a_0, \a_1, \ldots, \a_{n-1}, \gamma_1\}$ is a
basis of the root system of $\mathfrak{A}_1$ and $\{\b_1, \ldots,
\b_m, \gamma_2\}$ is a basis of the root system of
$\mathfrak{A}_2$.  Since a simple Lie algebra of type $D_r$, $r>4$
is not cuspidal by \cite{Fe}, then we have the following cases up
to a change of the basis of $\D$.

\medskip

\noindent{\it Case 1.} Let $\ga_i\in \Omega_v$, $i=1,2$. Then in
particular  $[\G_{\varepsilon},\G_{\varepsilon'}]v=0$ for all
$\varepsilon, \varepsilon'\in \dot{\Delta}^+_1$. It follows from
the proof of
 Lemma~\ref{le1}(i)  that in this case there
exists   a nonzero element $v'\in V$ which is  $S$-primitive for
$S=\pi\setminus \{\a_n\}$.

\medskip

\noindent{\it Case 2.} Let $m>1$. Suppose $\ga_1\in \Omega_v$,
$\ga_2\notin \Omega_v$ and $\b_2\in \Omega_v$.   Then in
particular all roots from
 $\dot{\D}_0^+$ containing $\b_2$ belong to $\Omega_v$. As in the proof
of Lemma~\ref{le1}(i) we can show that there exists   a nonzero
element $v'\in V$ such that $\Omega_v\subset \Omega_{v'}$ and
 $\Sigma_{\b}\subset \Omega_{v'}$ for all $\b\in \dot{\Delta}^+_1\setminus
 \{\a_n+\ldots + \a_{i}, \b_1+\a_n+ \ldots +\a_i, i=1, \ldots, n\}$.
 Hence $v'$  is  $S$-primitive for
$S=\pi\setminus \{\b_2\}$.

\medskip

\noindent{\it Case 3.} Let $m>1$. Suppose $\ga_1, \b_3\in
\Omega_v$ and  $\ga_2, \b_2\notin \Omega_v$. Then as in the
previous case there exists a nonzero element $v'\in V$ such that
$\Omega_v\subset \Omega_{v'}$ and
 $\Sigma_{\b}\subset \Omega_{v'}$ for all $\b\in \dot{\Delta}^+_1\setminus
 \{\a_n+\ldots + \a_{i}, \b_1+\a_n+ \ldots +\a_i, \b_2+\b_1+\a_n+ \ldots +\a_i,
 i=1, \ldots, n\}$.
 Hence $v'$  is  $S$-primitive for
$S=\pi\setminus \{\b_3\}$.

\medskip

\noindent{\it Case 4.}  Suppose $\ga_2\in \Omega_v$, $\ga_1\notin
\Omega_v$ and $\a_i\in \Omega_v$ for some $i=0, \ldots, n-1$. In
this case there exists a  $S$-primitive element with
$S=\pi\setminus \{\a_i\}$.

\medskip

\noindent{\it Case 5.} Let $m>1$.  Suppose $\ga_1, \ga_2\notin
\Omega_v$,  $\b_2\in \Omega_v$ and $\a_m\in \Omega_v$ for some
$m=0, \ldots, n-1$.   Then in particular all roots from
 $\dot{\D}_0^+$ containing $\b_2$ belong to $\Omega_v$.
  As in the proof of Lemma~\ref{le1} one can find a nonzero element
$v'\in X$ which is $(\ve, \pi)$-admissible for all $\ve\in
\dot{\D}_1^+$. Now consider the following  odd roots $z_{1i}=
\a_n+\ldots +\a_i$ and $z_{2i}= \b_1+\a_n+\ldots +\a_i$,  $i=1,
\ldots, n$. Define nonzero elements $v_{ji}$ by:
$v_{ji}=v_{j,i-1}$ if $X_{-z_{ji}}v_{j, i-1}=0$,
$v_{ji}=X_{-z_{ji}}v_{j,i-1}$ otherwise, $j=1,2$, $i=1, \ldots,
n$, $v_{11}=v'$, $v_{21}=v_{1n}$. Then $v_{2n}$ is a
$\P$-primitive element where $\P^0$ is spanned by $\pm \a_0,
\ldots, \pm\a_{m-1}, \pm\a_{m+1},\ldots,\pm\a_{n-1}, \pm\ga_1$, $\pm\b_1,
\pm\b_3, \ldots, \pm\b_m, \pm\ga_2$ and
$$\P^+=(\D_+\setminus (\P^0\cup \{z_{ji}, j=1,2, i=1, \ldots, n\}))\cup
\{-z_{ji}, j=1,2, i=1, \ldots, n\}.$$

\medskip

\noindent{\it Case 6.} Let $m>1$. Suppose $\ga_1, \ga_2,
\b_2\notin \Omega_v$,  $\b_3\in \Omega_v$ and $\a_m\in \Omega_v$
for some $m\in \{0, \ldots, n-1\}$. As in Case 5, there exists a
$\P$-primitive element, where $\P^0$ is spanned by $\pm \a_0,
\ldots, \pm\a_{m-1}, \pm \a_{m+1}, \ldots, \pm\a_{n-1}, \pm\ga_1$,
$\pm\b_1, \pm\b_2, \b_4, \ldots, \pm\b_m, \pm\ga_2$ and
$$\P^+=(\D_+\setminus (\P^0\cup \{d_{ji}, j=1,2,3, i=1, \ldots, n\}))\cup
\{-d_{ji}, j=1,2,3, i=1, \ldots, n\},$$ $d_{1i}= \a_n+\ldots
+\a_i$, $d_{2i}= \b_1+\a_n+\ldots +\a_i$,  $d_{3i}=
\b_2+\b_1+\a_n+\ldots +\a_i$, $i=1, \ldots, n$.

\medskip

\noindent{\it Case 7.} Suppose $\ga_1, \b_1\in \Omega_v$. Then
there exists a $S$-primitive element, $S=\pi\setminus \{\b_1\}$.

\medskip

\noindent{\it Case 8.} Suppose $\ga_1\notin \Omega_v$, $\b_1\in
\Omega_v$ and $\a_m\in \Omega_v$ for some $m\in \{0, \ldots,
n-1\}$. Then there exists a $\P$-primitive element, where $\P^0$
is spanned by $$\pm \a_0, \ldots, \pm\a_{m-1}, \pm\a_{m+1},
\ldots, \pm\a_{n-1}, \pm\ga_1, \pm\b_2, \ldots, \pm\b_m,
\pm\ga_2.$$

\medskip

\noindent{\it Case 9.} Let $m=1$, $\ga_1$ as before and
$\ga_2=-\b_1+\d$.  Changing the basis of $\D$ if necessary  we can
assume that we have one of the following subcases:

\medskip

\noindent{\it Case 9.1.} Let $\gamma_i\in \Omega_v$, $i=1,2$. Then
there exists a $S$-primitive element $v'\in V$ such that
$S=\pi\setminus \{\a_n\}$.

\medskip

\noindent{\it Case 9.2.} Let $\gamma_1\notin \Omega_v$ and
$\gamma_2, \a_i\in \Omega_v$
 for some $i\in \{0, \ldots, n-1\}$. As in
Case~4 there exists a  $S$-primitive element with $S=\pi\setminus
\{\a_i\}$.

This completes the proof of Proposition~\ref{pr4} for
$B(m,n)^{\hat{}}$.

\subsection{Case $B(0,n)^{\hat{}}$}
Let  $\G$ be of type $B(0,n)^{\hat{}}$, $\pi=\{\a_0, \ldots,
\a_n\}$, $\dot{\pi}_0=\{\a_1, \ldots, \a_{n-1}\}$,
$\dot{\pi}_1=\{\a_n\}$, $\dot{\D}^+_0=\{\a_i+\ldots + \a_j,
\a_i+\ldots + \a_j+2\a_{j+1}+\ldots +2\a_n, 1\leq i\leq j\leq
n-1\}$ $\cup \{2\a_{i}+\ldots +2\a_n, 1\leq i \leq n-1\}$,
 $\dot{\D}^+_1=\{\a_{i}+\ldots +\a_n, 1\leq i\leq n\}$.

\medskip

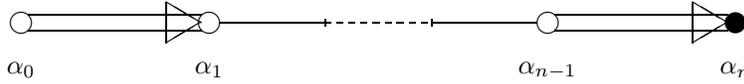
\unitlength=1.00mm
\special{em:linewidth 0.4pt}
\linethickness{0.4pt}
\begin{picture}(105.74,23.80)
\put(8.00,11.00){\circle{2.83}}
\put(33.00,11.00){\circle{2.83}}
\put(9.00,12.00){\line(1,0){23.00}}
\put(9.00,10.00){\line(1,0){23.00}}
\put(78.00,11.00){\circle{2.83}}
\put(103.00,11.00){\circle*{2.83}}
\put(34.40,11.00){\line(1,0){14.00}}
\put(79.00,12.00){\line(1,0){23.00}}
\put(79.00,10.00){\line(1,0){23.00}}
\put(62.60,11.00){\line(1,0){14.00}}
\put(48.40,11.00){\dashbox{0.8}(14.20,0.00)[cc]{}}
\put(8.00,4.67){\makebox(0,0)[cc]{$\alpha_0$}}
\put(33.00,4.67){\makebox(0,0)[cc]{$\alpha_1$}}
\put(78.00,4.67){\makebox(0,0)[cc]{$\alpha_{n-1}$}}
\put(103.00,4.67){\makebox(0,0)[cc]{$\alpha_n$}}
\put(32.00,11.00){\line(-5,3){4.67}}
\put(27.30,13.67){\line(0,-1){5.33}}
\put(27.30,8.33){\line(5,3){4.67}}
\put(102.00,11.00){\line(-5,3){4.67}}
\put(97.30,13.67){\line(0,-1){5.33}}
\put(97.30,8.33){\line(5,3){4.67}}
\end{picture}

\medskip

In this case $\cB_0\simeq C_n$.  Then $\{\a_0, \a_1, \ldots,
\a_{n-1}, 2\a_n\}$ is a basis of the root system of $\cB$.  Then
we have the following cases.

\medskip

\noindent{\it Case 1.}
Let $2\a_n\in \Omega_v$. Then in particular $\Omega_v$ contains all roots from
$\dot{\D}_0^+$ of the form $\b+2\a_n$ and
$[\G_{\varepsilon},\G_{\varepsilon'}]v=0$ for all $\varepsilon, \varepsilon'\in
\dot{\Delta}^+_1$. It follows from  the proof of
 Lemma~\ref{le1}(i)  that there
exists   a nonzero element $v'\in V$ which is  $S$-primitive for $S=\pi\setminus \{\a_n\}$.

\medskip

\noindent{\it Case 2.}
Suppose $2\a_n\notin \Omega_v$, $\a_0=-2\sum_{i=1}^{n}\a_i+\d\notin \Omega_v$ and
$\a_i\in \Omega_v$ for some $1\leq i\leq n-1$. Then in particular  all roots from
 $\dot{\D}_0^+$ containing $\a_i$ belong to $\Omega_v$. Same arguments as in the proof
of Lemma~\ref{le1}(i) show that there
exists   a nonzero element $v'\in V$ such that $\Omega_v\subset \Omega_{v'}$ and
 $\Sigma_{\b}\subset \Omega_{v'}$ for all $\b\in \dot{\Delta}^+_1\setminus \{\a_n, \a_{n}+\a_{n-1},
\ldots, \a_n+\ldots + \a_{i+1}\}$. Hence $v'$  is  $S$-primitive for $S=\pi\setminus \{\a_i\}$.

All other cases can be reduced to above by the change of a basis
of $\D$. This completes the proof of Proposition~\ref{pr4} for
$B(0,n)^{\hat{}}$.

\subsection{Case $D(m,n)^{\hat{}}$}
Let $\G$ be of type $D(m,n)^{\hat{}}$, $m\geq 2$,  $\pi=\{\a_0,
\ldots, \a_n,$  $\b_1, \ldots, \b_m\}$, $\dot{\pi}_0=\{\a_1,
\ldots, \a_{n-1}, \b_1, \ldots, \b_m\}$, $\dot{\pi}_1=\{\a_n\}$,
$\d=\a_0+2\sum_{i=1}^n +2\sum_{j=1}^{m-2} \b_j+\b_{m-1}+\b_m$.

\medskip

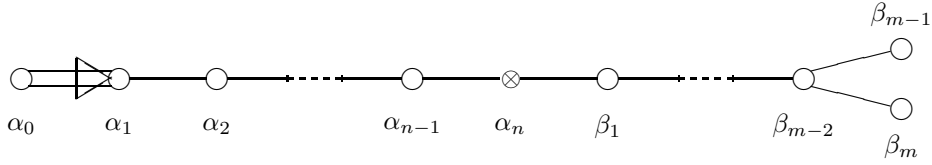
\unitlength=1.00mm
\special{em:linewidth 0.4pt}
\linethickness{0.4pt}
\begin{picture}(120.74,23.80)
\put(1.00,11.00){\circle{2.83}} \put(14.00,11.00){\circle{2.83}}
\put(2.0,12.00){\line(1,0){11.00}}
\put(2.00,10.00){\line(1,0){11.00}}
\put(27.00,11.00){\circle{2.83}} \put(53.00,11.00){\circle{2.83}}
\put(64.78,10.00){$\otimes$} \put(79.00,11.00){\circle{2.83}}
\put(105.00,11.00){\circle{2.83}} \put(118.00,7.00){\circle{2.83}}
\put(118.00,15.00){\circle{2.83}}
\put(15.40,11.00){\line(1,0){10.20}}
\put(54.40,11.00){\line(1,0){10.20}}
\put(67.40,11.00){\line(1,0){10.20}}
\put(28.40,11.00){\line(1,0){8.00}}
\put(43.60,11.00){\line(1,0){8.00}}
\put(88.40,11.00){\dashbox{0.8}(7.20,0.00)[cc]{}}
\put(80.40,11.00){\line(1,0){8.00}}
\put(95.60,11.00){\line(1,0){8.00}}
\put(36.40,11.00){\dashbox{0.8}(7.20,0.00)[cc]{}}
\put(1.00,4.67){\makebox(0,0)[cc]{$\alpha_0$}}
\put(14.00,4.67){\makebox(0,0)[cc]{$\alpha_1$}}
\put(27.00,4.67){\makebox(0,0)[cc]{$\alpha_2$}}
\put(53.00,4.67){\makebox(0,0)[cc]{$\alpha_{n-1}$}}
\put(66.00,4.67){\makebox(0,0)[cc]{$\alpha_n$}}
\put(79.00,4.67){\makebox(0,0)[cc]{$\beta_1$}}
\put(105.00,4.67){\makebox(0,0)[cc]{$\beta_{m-2}$}}
\put(118.00,19.50){\makebox(0,0)[cc]{$\beta_{m-1}$}}
\put(118.00,2.00){\makebox(0,0)[cc]{$\beta_m$}}
\put(13.00,11.00){\line(-5,3){4.67}}
\put(8.30,13.67){\line(0,-1){5.33}}
\put(8.30,8.33){\line(5,3){4.67}}
\put(106.00,12.00){\line(4,1){10.5}}
\put(106.00,10.00){\line(4,-1){10.5}}
\end{picture}

\medskip

Since $\cB_0\simeq C_n\oplus D_m$  we have that $\cB\simeq
\mathfrak{A}_1 + \mathfrak{A}_2$, where $\mathfrak{A}_1$ is an
affine Lie algebra of type $C_n^{(1)}$ while $\mathfrak{A}_2$ is
an affine Lie algebra of type $D_m^{(1)}$ for $m>3$, of type
$A_3^{(1)}$ for $m=3$ and of type $A_2^{(1)}$ for $m=2$. Denote
$\gamma_1=2\a_n+2\sum_{j=1}^{m-2} \b_j+\b_{m-1}+\b_m$ and
$\gamma_2=-\b_1-2\sum_{j=2}^{m-2} \b_j-\b_{m-1}-\b_m+\d$. Then
$\{\a_0, \a_1, \ldots, \a_{n-1}, \gamma_1\}$ is a basis of the
root system of $\mathfrak{A}_1$ and $\{\b_1, \ldots, \b_m,
\gamma_2\}$ is a basis of the root system of $\mathfrak{A}_2$. The
proof of Proposition~\ref{pr4} in this case is fully analogous to
the case of $B(m,n)^{\hat{}}$ and we leave the details out.

\subsection{Case $D(2,1;a)^{\hat{}}$}
Let $\G$ be of type $D(2,1;a)^{\hat{}}$, $\pi=\{\a_0, \ldots, \a_3\}$, $\dot{\pi}_0=\{\a_2, \a_3\}$,
$\dot{\pi}_1=\{\a_1\}$, $\dot{\D}^+_0=\{\a_2,\a_3, \a_2+2\a_1+\a_3\}$,
 $\dot{\D}^+_1=\{\a_1,\a_1 +\a_2, \a_1+\a_3, \a_1+\a_2+\a_3\}$.

\medskip

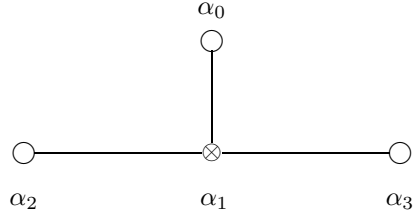
\unitlength=1.00mm
\special{em:linewidth 0.4pt}
\linethickness{0.4pt}
\begin{picture}(100.74,31.75)
\put( 25.00,11.00){\circle{2.83}}
\put( 48.60,10.20){$\otimes$}
\put( 75.00,11.00){\circle{2.83}}
\put( 50.00,26.00){\circle{2.83}}
\put( 25.00, 4.67){\makebox(0,0)[cc]{$\alpha_2$}}
\put( 50.33, 4.67){\makebox(0,0)[cc]{$\alpha_1$}}
\put( 75.00, 4.67){\makebox(0,0)[cc]{$\alpha_3$}}
\put( 50.00, 30.00){\makebox(0,0)[cc]{$\alpha_0$}}
\put( 26.40, 11.00){\line(1,0){22.20}}
\put( 51.40, 11.00){\line(1,0){22.20}}
\put( 50.00, 12.40){\line(0,1){12.20}}
\end{picture}

\medskip

Note that $\cB_0\simeq \mathfrak{sl}(2)\oplus
\mathfrak{sl}(2)\oplus \mathfrak{sl}(2)$. Consider the following
cases.

\medskip

\noindent{\it Case 1.}
Suppose that $\a_2, \a_3$ and  $\a_2+2\a_1+\a_3$ are not in $\Omega_v$.
Then by Lemma~\ref{le1}(i) there exists a nonzero
element $v'\in V$ which is $(\b, \pi)$-admissible for all
$\b\in \dot{\D}^+_1$.
Therefore $v'$ is $S$-primitive for $S=\pi\setminus \{\a_0\}$.

\medskip

\noindent{\it Case 2.}  Suppose now that $\a_2, \a_3, \a_0\notin \Omega_v$,
 where $\a_0=-\a_2-2\a_1-\a_3+\d$. Then by Lemma~\ref{le1}(i) there exists a nonzero
element $v'\in V$ which is
$\pi$-admissible and
$(\b, \pi)$-admissible for all
$\b\in \dot{\D}^+_1$. Denote the subspace of all such elements by $V'$. Consider
$w=X_{\a_1}v'$. Suppose $w\neq 0$. Then $w\in V'$ and  $X_{\a_1}w=0$. Next consider
$\a\in \dot{\D}^+_1$ of the least height such that $X_{\a}w\neq 0$ and repeat the argument above.
Note that $X_{\a}w\in V'$ and $X_{\b}X_{\a}w=0$ for any $\b\in \dot{\D}^+_1$. Eventually
we will find a nonzero $\tilde{v}\in V'$ such that $\dot{\D}^+_1\subset \Omega_{\tilde{v}}$, and hence
 $\tilde{v}$ is $S$-primitive for $S=\pi\setminus \{\a_1\}$.

All other cases can be reduced to Case 1 or Case 2 by a change of basis of $\D$. This completes
the proof of Proposition~\ref{pr4} for $D(2,1;a)^{\hat{}}$.

\subsection{Case $G(3)^{\hat{}}$}
Suppose that $\G$ is of type $G(3)^{\hat{}}$ with Dynkin diagram

\medskip

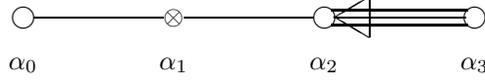
\unitlength=1.00mm
\special{em:linewidth 0.4pt}
\linethickness{0.4pt}
\begin{picture}(91.74,25.75)
\put( 20.00,11.00){\circle{2.83}}
\put( 38.60,10.00){$\otimes$}
\put( 60.00,11.00){\circle{2.83}}
\put( 80.00,11.00){\circle{2.83}}
\put( 20.00, 4.67){\makebox(0,0)[cc]{$\alpha_0$}}
\put( 40.00, 4.67){\makebox(0,0)[cc]{$\alpha_1$}}
\put( 60.00, 4.67){\makebox(0,0)[cc]{$\alpha_2$}}
\put( 80.00, 4.67){\makebox(0,0)[cc]{$\alpha_3$}}
\put( 61.40,11.00){\line(1,0){17.33}}
\put( 61.00,10.00){\line(1,0){18.00}}
\put( 61.00,12.00){\line(1,0){18.00}}
\put( 21.40,11.00){\line(1,0){17.20}}
\put( 41.40,11.00){\line(1,0){17.20}}
\put( 61.33,11.00){\line(5,3){4.67}}
\put( 66.00,13.87){\line(0,-1){5.53}}
\put( 66.00, 8.33){\line(-5,3){4.67}}
\end{picture}

\medskip

Let
$\pi=\{\a_0, \ldots, \a_3\}$, $\dot{\pi}_0=\{\a_2, \a_3\}$,
$\dot{\pi}_1=\{\a_1\}$, $\dot{\D}^+_0=\{\a_2,\a_3,\a_2+\a_3,2\a_2+\a_3,
3\a_2+\a_3, 3\a_2+2\a_3, 2\a_1+ 4\a_2+ 2\a_3\}$,
 $\dot{\D}^+_1=\{\a_1, \a_1 +\a_2, \a_1+ \a_2+\a_3, \a_1+2\a_2+\a_3, \a_1+3\a_2+\a_3,
\a_1+3\a_2+2\a_3, \a_1+4\a_2+2\a_3\}$. In this case $\cB_0\simeq
\mathfrak{sl}(2)\oplus {\bf G}_2$.

Denote $\b=\a_0+2\a_1+\a_2=-3\a_2-2\a_3+\d$. Then $\b, \a_2, \a_3$
is a basis of the root system of the affinization of ${\bf G}_2$.
Consider the following cases.

\medskip

\noindent{\it Case 1.}
Suppose that $\b \in \Omega_v$ and $-\a_0+\delta\in \Omega_v$. Then
$[\G_{\varepsilon},\G_{\varepsilon'}]v=0$ for all $\varepsilon, \varepsilon'\in
\dot{\Delta}^+_1$. It is not difficult to modify the proof of
 Lemma~\ref{le1}(i)  to show the
existence of  a nonzero
element $v'\in V$ which is  $S$-primitive for $S=\pi\setminus \{\a_1\}$.

\medskip

\noindent{\it Case 2.}  Suppose now that $\b\in \Omega_v$ but $-\a_0+\d\notin \Omega_v$.
Then in particular $v$ is $(\b, \pi)$-admissible
 for all $\b\in \dot{\D}^+_0$. Again using the same arguments as in the proof of
 Lemma~\ref{le1}(i)  we can find a nonzero
element $v'\in V$ which is $(\b, \pi)$-admissible
 for all $\b\in \dot{\D}\cap \D_+$. Hence $v'$ is
 $S$-primitive where $S=\pi\setminus \{\a_0\}$.

\medskip

\noindent{\it Case 3.}  Suppose that  $\a_2 \in \Omega_v$ and $-\a_0+\delta\in \Omega_v$.
 Denote $Y=\{\a_0, \a_3, \a_0+2\a_1+\a_2, \a_0+2\a_1+\a_2+\a_3\}$.
Then
 $\Omega_v$ contains all positive even roots except those in $Y$.
 Let $X$ be a subspace of $V^+$ of all  elements $w$ such that $(\D_0^+\setminus Y)\subset \Omega_w$.
Using the same algorithm as in the proof of Lemma~\ref{le1} we can find a nonzero element
$v_0\in X$ which is $(\ve, \pi)$-admissible for all $\ve\in \dot{\D}_1^+$. Finally, consider
the following  odd roots $\b_1=\a_1+4\a_2+2\a_3$, $\b_2=\a_1+3\a_2+2\a_3$,
$\b_3=\a_1+3\a_2+\a_3$, $\b_4=\a_1+2\a_2+\a_3$, $\b_5=-\a_1$, $\b_6=\a_1+\a_2$,
$\b_7=\a_1+\a_2+\a_3$, $\b_8=-\a_0-\a_1$ and define nonzero elements $v_i$ by:
$v_i=v_{i-1}$ if $X_{\b_i}v_{i-1}=0$ and $v_i=X_{\b_i}v_{i-1}$ otherwise, $i=1, \ldots, 8$.
Then $v_8\neq 0$ is a $\P$-primitive element where $\P^0=-Y\cup Y$ and
$$\P^+=(\D_+\setminus (Y\cup \{\a_1, \a_0+\a_1\}))\cup
\{-\a_1, -\a_0-\a_1\}.$$

\medskip

\noindent{\it Case 4.}  Suppose that  $\a_2 \in \Omega_v$ and $-\a_0+\delta\notin \Omega_v$.
Denote $Y=\{2\a_1+4\a_2+2\a_3, \a_1+2\a_2+\a_3,  \a_3, \a_0+2\a_1+\a_2,
 \a_0+2\a_1+\a_2+\a_3 \}$.
As in the previous case we find a nonzero $\P$-primitive element, where
 $\P^0=-Y\cup Y$ and $$\P^+=(\D_+\setminus (Y\cup \{\a_1, \a_1+\a_2, \a_1+\a_2+\a_3\}))\cup
\{-\a_1, -\a_1-\a_2, -\a_1-\a_2-\a_3\}.$$

\medskip

\noindent{\it Case 5.}  Suppose that  $\a_3 \in \Omega_v$ and $-\a_0+\delta\in \Omega_v$.
In this case there exists a nonzero $S$-primitive element, where $S=\pi\setminus \{\a_3\}$.

\medskip

\noindent{\it Case 6.}  Suppose that  $\a_3 \in \Omega_v$ and $-\a_0+\delta\notin \Omega_v$.
Denote $Y=\{2\a_1+4\a_2+2\a_3,  \a_2, \a_0+2\a_1+\a_2, \a_1+\a_2+\a_3, \a_1+2\a_2+\a_3,
\a_1+3\a_2+\a_3\}$.
In this case there exists a nonzero $\P$-primitive element, where $\P^0=-Y\cup Y$
and $\P^+=(\D_+\setminus (Y\cup \{\a_1, \a_1+\a_2\}))\cup
\{-\a_1, -\a_1-\a_2\}$.
This completes the proof of Proposition~\ref{pr4} for $G(3)^{\hat{}}$.

\subsection{Case $F(4)^{\hat{}}$}
Let $\G$ be of type $F(4)^{\hat{}}$ with Dynkin diagram

\medskip

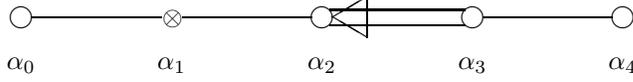
\unitlength=1.00mm
\special{em:linewidth 0.4pt}
\linethickness{0.4pt}
\begin{picture}(100.74,25.75)
\put( 10.00,11.00){\circle{2.83}}
\put( 28.77,10.00){$\otimes$}
\put( 50.00,11.00){\circle{2.83}}
\put( 70.00,11.00){\circle{2.83}}
\put( 90.00,11.00){\circle{2.83}}
\put( 10.00, 4.67){\makebox(0,0)[cc]{$\alpha_0$}}
\put( 30.00, 4.67){\makebox(0,0)[cc]{$\alpha_1$}}
\put( 50.00, 4.67){\makebox(0,0)[cc]{$\alpha_2$}}
\put( 70.00, 4.67){\makebox(0,0)[cc]{$\alpha_3$}}
\put( 90.00, 4.67){\makebox(0,0)[cc]{$\alpha_4$}}
\put( 71.40,11.00){\line(1,0){17.20}}
\put( 51.00,10.00){\line(1,0){18.00}}
\put( 51.00,12.00){\line(1,0){18.00}}
\put( 11.40,11.00){\line(1,0){17.20}}
\put( 31.40,11.00){\line(1,0){17.20}}
\put( 51.33,11.00){\line(5,3){4.67}}
\put( 56.00,13.87){\line(0,-1){5.53}}
\put( 56.00, 8.33){\line(-5,3){4.67}}
\end{picture}

\medskip

Let
$\pi=\{\a_0, \ldots, \a_4\}$, $\dot{\pi}_0=\{\a_2, \a_3, \a_4\}$,
$\dot{\pi}_1=\{\a_1\}$, $\dot{\D}^+_0=\{\a_2,\a_3, \a_4, \a_2+\a_3, 2\a_2+\a_3,
\a_3+\a_4, \a_2+\a_3+\a_4, 2\a_2+ \a_3+ \a_4, 2\a_2+2\a_3+\a_4, 2\a_1+3\a_2+2\a_3+\a_4\}$,
 $\dot{\D}^+_1=\{\a_1, \a_1 +\a_2, \a_1+ \a_2+\a_3, \a_1+2\a_2+\a_3, \a_1+\a_2+\a_3+\a_4,
\a_1+2\a_2+\a_3+\a_4, \a_1+2\a_2+2\a_3+\a_4, \a_1+3\a_2+2\a_3+\a_4\}$.
 In this case $\cB_0\simeq \mathfrak{sl}(2)\oplus B_3$.

  Let $\b=\a_0+2\a_1+\a_2$. Then
$\{\a_0, -\a_0+\delta\}\cup \{\a_2, \a_3, \b\}$ is a basis of $\cB$ and
 at least one element from each set is in $\Omega_v$. Hence
we have the following cases.

\noindent{\it Case 1.}
 If $-\a_0+\delta\in \Omega_v$ and $\b\in \Omega_v$ then
$[\G_{\varepsilon},\G_{\varepsilon'}]v=0$ for all $\varepsilon,
\varepsilon'\in \dot{\Delta}^+_1$. It follows from the proof of
 Lemma~\ref{le1}(i)  that there exists
  a nonzero
element $v'\in V$ which is  $S$-primitive for $S=\pi\setminus \{\a_1\}$.

\medskip

\noindent{\it Case 2.} If $\a_0\in \Omega_v$ and $\b\in \Omega_v$
then $v$ is $(\b, \pi)$-admissible
 for all $\b\in \dot{\D}^+_0$. The proof of
 Lemma~\ref{le1}(i) implies in this case that there exists a nonzero
element $v'\in V$ which  is
 $S$-primitive for $S=\pi\setminus \{\a_0\}$.
\medskip

\noindent{\it Case 3.} Let $-\a_0+\d\in \Omega_v$ and $\a_3\in
\Omega_v$. Again, as in the proof of
 Lemma~\ref{le1}(i) we can find  a nonzero
element $v'\in V$ which  is
 $S$-primitive for $S=\pi\setminus \{\a_3\}$.

\medskip

\noindent{\it Case 4.} Let $-\a_0+\d\in \Omega_v$ and $\a_2\in
\Omega_v$. Denote $Y=\{\a_3, \a_4, \a_3+\a_4, \b, \b+\a_3,
\b+\a_3+\a_4, \a_0\}$. In this case there exists a nonzero
$\P$-primitive element, where $\P^0=-Y\cup Y$ and
$\P^+=(\D_+\setminus (Y\cup \{\a_1,
-\a_1-3\a_2-2\a_3-\a_4+\d\}))\cup \{-\a_1,
 \a_1+3\a_2+2\a_3+\a_4-\d\}$.

\medskip

\noindent{\it Case 5.} Let $\a_0\in \Omega_v$ and $\a_3\in
\Omega_v$. Denote $Y=\{\a_2, \a_4, \a_1+\a_2+\a_3,
\a_1+\a_2+\a_3+\a_4, \a_1+2\a_2+\a_3,\a_1+2\a_2+\a_3+\a_4, -2\a_2-2\a_3-\a_4+\d,
2\a_1+3\a_2+2\a_3+\a_4\}$. In this case there exists a nonzero
$\P$-primitive element, where $\P^0=-Y\cup Y$ and
$\P^+=(\D_+\setminus (Y\cup \{\a_1, \a_1+\a_2\}))\cup \{-\a_1,
 -\a_1-\a_2\}$.

\medskip

\noindent{\it Case 6.} Let $\a_0\in \Omega_v$ and $\a_2\in
\Omega_v$. Denote $Y=\{\a_3, \a_4, \a_3+\a_4, \a_0+2\a_1+\a_2,
 \a_0+2\a_1+\a_2+\a_3, \a_0+2\a_1+\a_2+\a_3+\a_4, 2\a_1+3\a_2+2\a_3+\a_4\}$.
 In this case there
exists a nonzero $\P$-primitive element, where $\P^0=-Y\cup Y$ and
$\P^+=(\D_+\setminus (Y\cup \{\a_1, \a_1+\a_2, \a_1+\a_2+\a_3,
\a_1+\a_2+\a_3+\a_4\}))\cup \{-\a_1,
 -\a_1-\a_2, -\a_1-\a_2-\a_3, -\a_1-\a_2-\a_3-\a_4\}$.
This completes the proof of Proposition~\ref{pr4} for
 $F(4)^{\hat{}}$.

\section{Irreducible modules in $K(\G)$}

\noindent {\bf Proof of Theorem A}.
Let $V$ be an irreducible module in $K(\G)$. Then $V$ contains a $\P$-primitive
element for some  parabolic subset $\P\subset \D$ by Proposition~\ref{pr4}.
 Therefore
 $V\simeq L_{\P}(N')$ by Proposition~\ref{pr2} for  some irreducible
 module $N'$ over a finite-dimensional Lie subsuperalgebra $\G_{\P}^0$. By \cite{DMP},Theorem 6.1,
module $N$ is parabolically induced from a cuspidal module, that is there exists a
 parabolic subset $\P'$ of the root system of $\B=\G_{\P}^0$
 and an irreducible cuspidal $\B$-module $N$ such that $N'$ is a quotient of
${\rm ind}(\B_{\P'}^0\oplus \B_{\P'}^+, \B; N)$.  Combining the parabolic subsets $\P$
and $\P'$ we obtain a parabolic subset $\tilde{\P}$ of $\D$
 such that $V\simeq L_{\tilde{\P}}(N)$, a unique irreducible quotient of
${\rm ind}(\G_{\tilde{\P}}, \G; N)$.
 This completes the proof of Theorem~A.

Theorem A reduces the classification of irreducible modules in $K(\G)$ to the classification
of irreducible cuspidal modules over cuspidal Lie subsuperalgebras. 

\medskip

\noindent {\bf Proof of Theorem B}.
Let $V$ be an irreducible module in $K(\G)$. Then $V$ is isomorphic to
$L_{\P}(M)$ by Theorem A, where $\P$ is a parabolic subset of $\D$
and $M$ is an irreducible $\G^0_{\P}$-module.

 If $\G$ is of type
$A(m,n)^{\hat{}}$ then
 it follows from the proof of Proposition~\ref{pr4} that  $V\simeq
 L^{\pi_0}(M)$, where
$M$ is  irreducible  $\g_0$-module. Since $\g_0$ is
finite-dimensional reductive Lie algebra,
 then $M\simeq L^S_{\g_0}(N)$ for some subset $S\subset \pi_0$ and some
irreducible cuspidal $\g(S)$-module $M$ by \cite{Fe}, implying
(1).

Suppose that $\G$ is of type $C(n)^{\hat{}}$. It follows from the
proof of Proposition~\ref{pr4} that $V$ contains a $S$-primitive
element of weight $\l$ where either $S=\pi_0$ or $S=\pi\setminus
\{\a_k\}$ for some $2\leq k\leq n-1$. In the first case $V\simeq
L^{\pi_0}(M)$, where $M$ is irreducible weight module over
$\g_0\simeq \mathfrak{sp}(2n-2)\oplus \C$. Hence, by \cite{Fe},
there exists $S\subset \pi_0$ and an irreducible cuspidal
$\g(S)$-module $N$ such that $M\simeq L^S_{\g_0}(N)$ and $V\simeq
L^S(N)$. Now consider the case $S=\pi\setminus \{\a_k\}$. By
Proposition~\ref{pr1}, $V\simeq L^S(M)$ where $M=\sum_{\nu\in
Q_S}V_{\nu+\l}$  is an irreducible module over a
finite-dimensional Lie subsuperalgebra $\G(S)$. The superalgebra
$\G(S)$ is a direct sum of
 $G_1$ and $G_2$ where $G_1\simeq \mathfrak{osp}(2,2k-2)$ is a Lie superalgebra generated
by $\G_{\b}$, $\b\in \{\a_0, \a_1, \ldots, \a_{k-1}\}$ and
$G_2\simeq \mathfrak{sp}(2n-2k-2)$ is a Lie algebra generated by
$\G_{\b}$, $\b\in \{\a_{k+1}, \ldots, \a_n\}$. Then $M$ is
isomorphic to a tensor product $N_1\otimes N_2$ where $N_i$ is
irreducible $G_i$-module, $i=1,2$. Module $N_2$ is  parabolically
induced from  a cuspidal module  over a subalgebra of $G_2$.  Let
$\tilde{\pi}$ be a basis of the root system of $G_1$. Then $N_1$
contains a $\tilde{\pi}_0$-primitive element by
Lemma~\ref{le1}(ii), and therefore $N_1$ is  parabolically induced
from a cuspidal module over a  subalgebra of the even part of
$G_1$, implying (2).

Let $\G$ be of type $B(m,n)^{\hat{}}$ ($m>0$). By
Proposition~\ref{pr4}, one can choose $\P$ in such a way that
$\G^0_{\P}$ is isomorphic to one of the (super)algebras of the
following types: $C_n\oplus B_m$, $D(2,n)\oplus B_{m-2}$,
$D(3,n)\oplus B_{m-3}$, $C_k\oplus B(m,n-k), k=0, \ldots, n-1$,
$C_i\oplus C_{n-i}\oplus A_1\oplus A_1\oplus B_{m-2}$, $C_i\oplus
C_{n-i}\oplus A_2\oplus B_{m-3}$, $C(n+1)\oplus B_m$ and
$C_i\oplus C(n-i)\oplus B_m$. The only simple cuspidal subalgebras
of $B_m$ are those of type $A$ and $C_2$ by \cite{Fe}. The
superalgebras $D(2,n)$, $D(3,n)$ and $B(m,n-k)$ with $m=1,2$ are
cuspidal by \cite{DMP},7.3. Note that $\mathfrak{C}=B(m,n-k)$ is
not cuspidal if $m>2$. The reason is that in this case
$\mathfrak{C}$ contains a subalgebra of type $B_m$ which is not
cuspidal, and a subgroup generated by the  even roots of
$\mathfrak{C}$ has a finite index in the root lattice. Also note
that $C(k)$ is not cuspidal by \cite{DMP1},7.3 ( also by (2)).
Hence the statement (3) follows.

Suppose that $\G$ is of type $B(0,n)^{\hat{}}$. Then
 $V$ contains a $S$-primitive element,  where $S=\pi\setminus
\{\a\}$, $\a\in \pi$, by the proof of Proposition~\ref{pr4}. If
$\a\in \pi_0$ then $\G(S)\simeq \mathfrak{sp}(2i)\oplus
\mathfrak{osp}(1,2(n-i))$, where $\a=\a_i$, $i\in \{0, \ldots,
n-1\}$. If $\a$ is an odd root then $\G(S)\simeq
\mathfrak{sp}(2n)$. Since  $\mathfrak{osp}(1,2(n-i))$ is cuspidal
by \cite{DMP},7.3, and any semisimple subalgebra of
$\mathfrak{sp}(2n)$ is cuspidal by \cite{Fe}, then the statement
(4) follows.

The proof of (5) is similar to the proof of (3). Just note that
$\mathfrak{osp}(2m,2l)$ is cuspidal only if $m=2$ or $m=3$ and
$so(2k)$ is cuspidal only when $k\leq 3$.

The statement (6) is immediate.

If $\G$ is of type $G(3)^{\hat{}}$ then  it follows from the proof
of Proposition~\ref{pr4} that $\P$ can be chosen in such a way
that $\G^0_{\P}$ is isomorphic to one of the following:
$$ \mathfrak{sl}(2)\oplus
 {\bf G}_2,  \mathfrak{sl}(2)\oplus
\mathfrak{sl}(3), \mathfrak{sl}(3)\oplus \mathfrak{osp}(1,2),
 \mathfrak{osp}(4,2),  \mathfrak{sl}(2)\oplus \mathfrak{osp}(3,2), G(3).$$

 Simple cuspidal subsuperalgebras of $G(3)$ are isomorphic to $\mathfrak{sl}(2),  \mathfrak{osp}(1,2)$
and  $\mathfrak{osp}(3,2)$ by \cite{DMP},7.7 and \cite{DMP1},
while any cuspidal subalgebra of ${\bf G}_2$ is isomorphic to
$\mathfrak{sl}(2)$ by \cite{Fe}. Since  $\mathfrak{osp}(4,2)$ is
cuspidal
  by \cite{DMP},7.8,
the statement (7) follows.

 Let $\G$ be of type $F(4)^{\hat{}}$. It follows from the proof of Proposition~\ref{pr4}
 that $\P$ can be
chosen such that $\G^0_{\P}$ is isomorphic to one of the
following:
$$\mathfrak{sl}(2)\oplus
 \mathfrak{so}(7), \mathfrak{sl}(2)\oplus
\mathfrak{osp}(4,2), \mathfrak{sl}(2)\oplus
\mathfrak{sl}(4), F(4).$$
 Hence the statement (8) follows from
 \cite{Fe} and \cite{DMP},7.8.
This completes the proof of Theorem~B.

 From Theorem~B we immediately obtain a classification of all cuspidal Levi subsuperalgebras.

\begin{corol}\label{cor4}
Any simple cuspidal Levi subsuperalgebra of affine Lie superalgebra is isomorphic to
 one of the following:
\begin{itemize}
\item[(1)] Simple Lie subalgebra  of type $A$ or $C$;
\item[(2)] $\mathfrak{osp}(m,2n)$, $m=1,3,4,5,6$;
\item[(3)] $D(2,1;a)$.
\end{itemize}
\end{corol}

If $M$ is a module over a Lie algebra $\mathfrak{A}$ then we
denote by $inj M$ the set of all roots for which the corresponding
root element acts injectively on $V$.

\begin{corol}\label{cor5}
Let $V\simeq L_{\P}(N)$ be an irreducible module in $K(\G)$, where $\P$ is a parabolic subset of $\D$
and $N$ is a cuspidal $\G^0_{\P}$-module. Then either $\P^0=\emptyset$ or
there exists an even root $\a\in \D_0$ such that $\a\in inj N$, that is there exists a root element
$x\in \G_{\a}$ which acts injectively on $N$.
\end{corol}

\begin{proof}
Suppose $\P^0\neq\emptyset$. If $\P^0$ contains a root $\a\in \dot{\pi}_0$ then $\a\in inj V$ by definition.
 Let $\P^0\cap \dot{\pi}_0=\emptyset$. Then  $\P^0$ contains an odd root $\b$ and $2\b\in \D_0$. If
$2\b\notin inj N$ then $N$ contains a nonzero  element
$v$ such that $\G_{2\b}v=0$. Then $v$ is a $\tilde{\P}$-primitive
element with $\P^+\varsubsetneq\tilde{\P}^+$ which is a
contradiction.
\end{proof}

\begin{remark}
\begin{itemize}
\item[(1)] Combining Theorem~B with Mathieu's classification of cuspidal modules over simple
Lie algebras of type $A$ and $C$ we obtain a complete classification of irreducible modules in $K(\G)$
when $\G$ is of type $A(m,n)^{\hat{}}$ and $C(n)^{\hat{}}$.
 Note that our reduction in these cases   is  fully independent of the results of \cite{DMP}.
\item[(2)]   If $\G$ is not of type $A(m,n)^{\hat{}}$ and $C(n)^{\hat{}}$ then
 our reduction relies on the classification of cuspidal Levi subsuperalgebras
in \cite{DMP}. In this case the classification of irreducible modules in $K(\G)$ is reduced to a
finite indeterminacy by \cite{DMP},Proposition 6.3.
\item[(3)] In cases  $C(n)^{\hat{}}$ and $G(3)^{\hat{}}$ the irreducible modules in $K(\G)$
are parabolically induced but
not always   isomorphic to modules of type $L^S(N)$.
\item[(4)] Classification of non-cuspidal irreducible modules in $K(\G)$ for $D(2,1;a)$ is based
on the classification of cuspidal $\mathfrak{sl}(2)$-modules and is, therefore, complete.
\item[(5)] Classification of all irreducible modules over $\mathfrak{osp}(1,2)$ was obtained in
 \cite{BO}.
\item[(6)] Classification of  irreducible modules in $K(\G)$ will
be complete after the classification of cuspidal modules over
$\mathfrak{osp}(m,2n)$, $m=1,3,4,5,6$ and $D(2,1;a)$. \item[(7)]
Cuspidal Lie subsuperalgebras described in Corollary~\ref{cor4}
coincide with cuspidal superalgebras  in the finite-dimensional
setting which were classified in \cite{DMP1}.
\end{itemize}
\end{remark}

\section{Integrable modules over affine Kac-Moody algebras}

Let $\G$ be affine Lie algebra, $\Delta$ its root system,
$\g$ the underlying finite-dimensional simple Lie
subalgebra of $\G$, $\dot{\Delta}$ the root system of $\g$.
A weight $\G$-module is {\it integrable} if all $\G_{\a}$,
$\a\in \D^{\re}$, are locally finite on $V$.
We will now recover the classification of
 irreducible integrable modules in $K(\G)$ which is due to Chari.

\begin{prop}\label{pr10}
\begin{itemize}
\item[(i)] If $\g$ is locally finite on an irreducible module $V$ in $K(\G)$ then $V$ is
a highest weight module.
\item[(ii)] (\cite{C}) Any  irreducible integrable module in $K(\G)$ is a highest weight module.
\end{itemize}
\end{prop}

\begin{proof}
Let $V$ be an irreducible module in $K(\G)$ on which $\g$ is locally finite. By Theorem~\ref{thm1}(ii),
$V\simeq L^S_{\G}(N)$ for some basis $\pi$ of $\Delta$, a proper subset $S\subset \pi$ and an irreducible cuspidal
$\G(S)$-module $N$. Suppose  $S\neq \emptyset$, hence $V$ is not a highest weight module. Since $\g$ is locally finite then $S\cap  \dot{\pi}=\emptyset$, where $\dot{\pi}$ is a basis of $\dot{\Delta}$. 
Suppose that $S$ contains a root $\alpha-m\delta$ where $\alpha\in \dot{\Delta}$, $m\in \Z$.   
Let $\cB(\a-m\delta)\simeq \mathfrak{sl}(2)$ a subalgebra of $\G$ generated
by $\G_{\pm (\a-m\delta)}$. Fix a nonzero weight vector $v\in N$. Then $U(\cB(\a-m\delta))v\subset N$ 
is a torsion free $\cB(\a-m\delta)$-module. Without loss of generality we can assume that $\G_{\a}v=0$ and $\G_{k\delta}v=0$ for all $k>0$. Under such assumption $m$ must be positive. Consider now a subalgebra $\cB(\a)\simeq \mathfrak{sl}(2)$  of $\G$ generated
by $\G_{\pm \a}$. 
Let $V_1=U(\cB(\a))v$. Then $V_1$ is a finite-dimensional $\mathfrak{sl}(2)$-module. Assume $\dim V_1=r$. Choose a basis  element $X_{\a+k\delta}$ in each space  $\G_{\pm (\a+k\delta)}$, $k\in \Z$ and 
consider $w_1=X_{-\a+m\delta}v$. Let   $V_2=U(\cB(\a))w_1$. Note that $\G_{k\delta}w_1=\G_{\alpha}w_1=0$ for all $k>0$ due to the induced structure of $V$. 
 Assume now that $\G$ is not of type $A_{2n}^{(2)}$. Then
 $[X_{-\a}, X_{-\a+m\delta}]=0$ and hence
$X_{-\a}^r w_1=0.$ Therefore, $V_2$ is a finite-dimensional module of dimension $\leq r$. Applying the same arguments to nonzero elements $w_i=X_{-\a+m\delta}^i v$, $i=1,\ldots, r+1$, we construct $r+1$ nonzero modules of dimension $\leq r$ over $\mathfrak{sl}(2)$. Since all these module have different highest weights as $\mathfrak{sl}(2)$-modules we obtain a contradiction. 

Now consider the case when $\G$ is  of type $A_{2n}^{(2)}$. In this case we have 
$$[X_{-\a},[X_{-\a}, X_{-\a+m\delta}]]=0.$$ Hence $X_{-\a}^{r+1} w_1=0.$ Therefore, $V_2$ is a finite-dimensional module of dimension $\leq r+1$. But the difference 
of highest weights of $V_1$ and $V_2$ is $2$, thus $V_2$ can not have dimension $r+1$. It has dimension $\leq r$ and the same arguments as above lead to a contradiction. Thus (i) follows. Statement (ii) is an immediate corollary of (i).
\end{proof}

\section{Weakly integrable modules in $K(\G)$}

Let $\G=\G_0\oplus \G_1$ be a non-twisted affine Lie superalgebra and $\g=\g_0\oplus \g_1$
the underlying basic classical Lie superalgebra. Let also $\g_0=\sum_{j=0}^s g_{0j}$
where $\g_{00}$ is abelian and $\g_{0i}$, $i=1,\ldots, s$, are simple Lie algebras ($s\leq 3$).
  Following \cite{KW2} we define {\it weakly integrable} $\G$-module $V$ 
as a module which is integrable over the affinization $\G_{0j}$ of
$\g_{0j}$ for some $j\geq 1$, and on which $\g_0$ is locally finite. Of course, for
$C(n)^{\hat{}}$ a weak integrability coincides with the usual   integrability since $s=1$ in this case.
Theorem~C states that the only weakly integrable irreducible modules in $K(\G)$
are highest weight modules which is a generalization of Proposition~\ref{pr10} for affine Lie superalgebras.

\begin{theorem}\label{theorem11}
Let $V$ be an irreducible module in $K(\G)$ such that $\g_0$ is locally finite on $V$. Then $V$ is
a highest weight module.
\end{theorem}

\begin{proof}
 By Theorem A there exists a parabolic subset $\P\subset \D$
 and an irreducible cuspidal $\G^0_{\P}$-module $N$ such that $V\simeq L_{\P}(N)$.
 Since $\g_0$ is locally finite on $V$ then
  $\G^0_{\P}\neq \G$. Suppose first that $\P^0\neq \emptyset$.
 Since $N$ is cuspidal then by Corollary~\ref{cor5} there exists an even root $\a$ such that
 $N$ is cuspidal over a subalgebra
 $\cB(\a)\simeq \mathfrak{sl}(2)$,
 generated by $\G_{\pm \a}$.
If $\a\in \dot{\D}_0$ then $\g_0$ is not locally finite on $N$, and hence on $V$,
 which is a
 contradiction. Suppose now that $\a=\phi-m\delta$ for some
$\phi\in \dot{\D}_0$ and $m\in \Z$. Without loss of generality we may assume that
$\phi\in \P^+$ (changing $\a$ to $-\a$ if necessary) and $\delta\in \P^+$. It implies $m >0$.
 As in the proof of Proposition~\ref{pr10} we can consider  a subalgebra
 $\cB(\phi- m\delta)\simeq \mathfrak{sl}(2)$  generated by $\G_{\pm (\phi - m\delta)}$. 
Let $M=U(\cB(\phi- m\delta))v$ for some nonzero weight element $v\in N$
and
$\tilde{M}=U(\g_0)M\subset V$.
Then $\tilde{M}$ is a  $\G(\phi)$-module, where $\G(\phi)\simeq \mathfrak{sl}(2)$ is generated by 
$\G_{\pm \phi}$.  
Applying the same argument as in the proof of Proposition~\ref{pr10} we conclude that
the action of  $\G_{-\phi}=\g_{-\phi}$ is not locally finite on
$\tilde{M}$ which is  a contradiction.
Therefore we can assume that $\P^0= \emptyset$. Let $v$ be a nonzero
$\P$-primitive weight element of $V$. Fix $j$ and let $\G_{0j}$ be the
affinization of $\g_{0j}$.
Let $\D_{0j}$ be the root system of $\G_{0j}$ and $\D_{0j}^+=\D_{0j}\cap \D_+$. 
Consider a $\G_{0j}+\H$-submodule $V'=U(\G_{0j})v$ of $V$. Since
 $v$ is a $\P$-primitive element then $V'$ is a homomorphic image of the Verma type module
${\rm ind}(\H\oplus (\G_{0j})^+_{\P},\G_{0j}; \C)$ with a nonzero central charge, 
where $(\G_{0j})^+_{\P}=\G_{0j}\cap \G_{\P}^+$. 
Suppose that $\delta\in \P$ and $w\P\cap \D_{0j}^+\neq \D_{0j}^+$ for any $w$ in the Weyl group of $\Delta_{0j}$.
Then $V'$  has some infinite-dimensional
weight spaces by \cite{F1},Theorem 5.14. But this is a contradiction since $V'$ is an object of
the category $K(\G)$. Therefore, there exists an element $w$ of the Weyl group of $\G$ such that 
$w\P\cap \D_{0j}^+= \D_{0j}^+$ for all $j$. It follows immediately that
$\Sigma_{\b}\subset \Omega_v$ for all $\b\in \dot{\D}_0^+(w^{-1}\pi)$.  Applying Lemma~\ref{le1},(iii), we obtain that $V$
 is a highest weight module.
 This completes the proof.
\end{proof}

Theorem C follows immediately from Theorem~\ref{theorem11}.

\begin{remark}
\begin{itemize}
\item[1.]
 Together with the Kac-Wakimoto classification of  integrable irreducible highest
weight modules Theorem C  gives a complete classification of
weakly integrable irreducible modules in $K(\G)$.

\item[2.] Using \cite{DMP},Theorem~6.1 one can show that
Theorem~\ref{theorem11} also holds for basic classical Lie
superalgebras. Therefore any irreducible integrable $\g$-module
with finite-dimensional weight spaces is highest weight.

\item[3.] If we relax the requirement of local finiteness of
$\g_0$ in the definition of weakly integrability then we get
non-highest weight "integrable" modules in $K(\G)$. Such modules
are partially integrable in the sense of Dimitrov and Penkov
(\cite{DP}).

\end{itemize}
\end{remark}

\section{Conjectures}

In  conclusion we formulate several conjectures 
for affine Lie algebras.

 Let $\cB$ be an affine Lie algebra with a Cartan
subalgebra $H$. Let $\mathfrak{k}(\cB)$ be the category of all
weight $\cB$-modules, $\tilde{K}(\cB)$ the full subcategory of
$\mathfrak{k}(\cB)$ consisting of modules
 with finite-dimensional weight subspaces and
 $K_0(\cB)$  the full subcategory of $\tilde{K}(\cB)$ of modules of zero level.

 A $\cB$-module $V\in \mathfrak{k}(\cB)$ is {\em dense} if the set of weights is
a  coset of $H^*/Q$.
 The following conjecture was formulated in \cite{F1}:

\bigskip

\noindent{\bf Conjecture 1.}  Irreducible module $V\in
\mathfrak{k}(\cB)$ is dense if and only if it is cuspidal.

\bigskip

This conjecture would reduce the classification problem to the
classification of irreducible weight modules with
infinite-dimensional weight spaces
 over the Heisenberg subalgebra
of $\cB$  and irreducible dense modules over affine Lie
subalgebras of $\cB$ (cf. \cite{F1}). Irreducible  modules with
infinite-dimensional weight spaces
 over the Heisenberg algebra  necessarily have a nonzero level. Some examples were
constructed in \cite{BBF}. Families of  dense $\cB$-modules with
infinite-dimensional weight spaces were constructed in \cite{CP}.
These classification problems are still open.

 Theorem~\ref{thm1}(ii) confirms Conjecture~$1$ for the category $K(\cB)$ and describes
 irreducible
 modules there. On the other hand all irreducible   modules
 in $K_0(\cB)$ which are not irreducible over $[\cB,\cB]$ were classified in
 \cite{R}.
 Such modules are called {\em Loop} modules.
  We believe that  the following holds:

\bigskip

\noindent{\bf Conjecture 2.} An irreducible module $V\in K_0(\cB)$
is either a Loop module or is parabolically induced from a
standard parabolic subalgebra of $\cB$ (cf. Section 4).

\bigskip

Conjecture~2 together with  Theorem~\ref{thm1}(ii) would complete
a classification of irreducible modules in $\tilde{K}(\cB)$.

It was shown in \cite{F} that Conjecture~$1$ holds for modules
over $\mathfrak{sl}(2)^{\hat{}}$ with at least one
finite-dimensional weight space. This suggests the following
weaker version of Conjecture~1:

\bigskip

\noindent{\bf Conjecture 3.} Irreducible module $V\in
\mathfrak{k}(\cB)$ with at least one finite-dimensional weight
space is dense if and only if it is cuspidal.

\bigskip

Irreducible integrable modules in $\mathfrak{k}(\cB)$ are
quotients of  tensor product of a standard highest weight and a
dual to standard modules by \cite{CP}.

\bigskip

\noindent{\bf Conjecture 4.} Irreducible integrable module $V\in
\mathfrak{k}(\cB)$ of  a nonzero level is  highest weight.


\bigskip




\begin{thebibliography}{99}
\bibitem{BO} {\em V.Bavula, F. van Oystaeyen}, The simple modules of the Lie superalgebra
$\mathfrak{osp}(1,2)$, J. Pure and Applied Algebra. 150 (2000), 41-52.
\bibitem{BBF} {\em V.Bekkert, G.Benkart, V.Futorny}, Weyl algebra
modules,  Kac-Moody Lie algebras and related topics,  17--42, Contemp. Math., 343, Amer. Math. Soc., Providence, RI, 2004.
\bibitem{C} {\em V.Chari}, Integrable representations of affine Lie algebras,
Invent. Math. 85 (1986), no.2, 317-335.
\bibitem{CP} {\em V.Chari, A.Pressley}, New unitary
representations of loop groups, Math. Ann. 275 (1986), 87-104.
\bibitem{DMP} {\em I.Dimitrov, O.Mathieu, I.Penkov}, On the structure of weight modules,
Trans. Amer. Math. Soc. 352 (2000), 2857-2869.
\bibitem{DMP1} {\em I.Dimitrov, O.Mathieu, I.Penkov},  Errata,
Trans. Amer. Math. Soc. (to appear).
\bibitem{DP} {\em I.Dimitrov, I.Penkov}, Partially and fully integrable modules over Lie superalgebras,
Studies in Advanced Mathematics (series editor S.-T. Yau), 4, AMS and Internatl. Press 1997, 49-67.
\bibitem{Fe}
{\em S.Fernando}, Lie algebra modules with finite dimensional
 weight spaces I, Trans. Amer.  Math. Soc. 322 (1990), 757-781.
\bibitem{F} {\em V.Futorny}, Irreducible non-dense
$A_1^{(1)}$-modules, Pacific J. of Math. 172 (1996), 83-99.
\bibitem{F1} {\em V.Futorny}, Representations of Affine Lie algebras,
Queen's Papers in Pure and Applied Mathematics, v. 106, 1997, Kingston, Canada.
\bibitem{FT}
{\em V.Futorny, A.Tsylke}, Classification of irreducible nonzero level modules
with finite-dimensional weight spaces for affine Lie algebras, J. Algebra 238 (2001), 426-441.
\bibitem{K} {\em V.Kac}, Lie superalgebras, Adv. Math. 26 (1977), 8-96.
\bibitem{KW1} {\em V.Kac, M.Wakimoto}, Integrable highest weight modules over affine Lie superalgebras
and number theory, Lie Theory and Geometry, Birkhauser, Progress in Mathematics 123 (1994), 415-456.
\bibitem{KW2} {\em V.Kac, M.Wakimoto}, Integrable highest weight modules over affine Lie superalgebras
and Appell's function, Comm.Math.Phys. 215 (2001), 631-682.
\bibitem{M} {\em O.Mathieu}, Classification of irreducible weight modules, Ann. Inst. Fourier 50 (2000), 537-592.
\bibitem{R}{\em S.Rao}, Classification of Loop modules with
finite-dimensional weight spaces, Math. Ann. 305 (1996), 651-663.
\bibitem{RZ}{\em S.Rao, K.Zhao}, On integrable representations for toroidal Lie superalgebras,
Contemporary Mathematics 343  (2004), 243-261.
\end{thebibliography}
\end{document}